%% file: lanczos-paper.tex
\documentclass[12pt]{article}

\pdfoutput=1

\usepackage{amsmath}
\usepackage{graphicx,psfrag,epsf}
\usepackage{enumerate}
\usepackage{array}
\usepackage{natbib}
\usepackage{url} 

\usepackage{lipsum}
\usepackage{amsfonts}
\usepackage{graphicx}
\usepackage{csquotes}
\usepackage{epstopdf}
\usepackage{algorithmic}
\ifpdf
  \DeclareGraphicsExtensions{.eps,.pdf,.png,.jpg}
\else
  \DeclareGraphicsExtensions{.eps}
\fi

\input{pauls-macros}

\usepackage{amsmath,mathtools}
\usepackage{tabulary,booktabs}
\usepackage{thmtools,thm-restate}
\usepackage{algorithm}
\usepackage{subfig}
\usepackage{bm}
\usepackage{color}
\usepackage{hyperref}
\usepackage{enumitem}
\setlist[enumerate]{itemsep=0mm}
\setlist[itemize]{itemsep=0mm}

\newtheorem{assume}{Assumption}

\newtheorem{remark}{Remark}

\newcommand{\blind}{1}

\begin{document}

\def\spacingset#1{\renewcommand{\baselinestretch}%
{#1}\small\normalsize} \spacingset{1}

\if1\blind
{
  \title{\bf Gauss-Christoffel quadrature for inverse regression: applications to computer experiments}
  \author{Andrew Glaws%
  \thanks{Department of Computer Science, University of Colorado, Boulder, CO (\url{andrew.glaws@colorado.edu}).}
    \and 
    Paul G. Constantine%
  \thanks{Department of Computer Science, University of Colorado, Boulder, CO (\url{paul.constantine@colorado.edu}).}}
  \date{}
  \maketitle
}\fi

\if0\blind
{
  \bigskip
  \bigskip
  \bigskip
  \begin{center}
    {\LARGE\bf Gauss-Christoffel quadrature for inverse regression: applications to computer experiments}
\end{center}
  \medskip
}\fi

\bigskip
\begin{abstract}
Sufficient dimension reduction (SDR) provides a framework for reducing the predictor space dimension in regression problems. We consider SDR in the context of deterministic functions of several variables such as those arising in computer experiments. In this context, SDR serves as a methodology for uncovering ridge structure in functions, and two primary algorithms for SDR---sliced inverse regression (SIR) and sliced average variance estimation (SAVE)---approximate matrices of integrals using a sliced mapping of the response. We interpret this sliced approach as a Riemann sum approximation of the particular integrals arising in each algorithm. We employ well-known tools from numerical analysis---namely, multivariate tensor product Gauss-Christoffel quadrature and orthogonal polynomials---to produce new algorithms that improve upon the Riemann sum-based numerical integration in SIR and SAVE. We call the new algorithms \emph{Lanczos-Stieltjes inverse regression} (LSIR) and \emph{Lanczos-Stieltjes average variance estimation} (LSAVE) due to their connection with Stieltjes' method---and Lanczos' related discretization---for generating a sequence of polynomials that are orthogonal to a given measure. We show that the quadrature-based approach approximates the desired integrals, and we study the behavior of LSIR and LSAVE with three numerical examples. As expected in high order numerical integration, the quadrature-based LSIR and LSAVE exhibit exponential convergence in the integral approximations compared to the first order convergence of the classical SIR and SAVE. The disadvantage of LSIR and LSAVE is that the underlying tensor product quadrature suffers from the curse of dimensionality---that is, the number of quadrature nodes grows exponentially with the input space dimension. Therefore, the proposed approach is most appropriate for deterministic functions with fewer than ten independent inputs.
\end{abstract}

\noindent%
\vfill
\hfill {\tiny technometrics tex template (do not remove)}

\newpage
\spacingset{1.45}

\input{sec1-intro}
\input{sec2-IR_methods}
\input{sec3-comp_exp}
\input{sec4-num_tools}
\input{sec5-LS_methods}
\input{sec6-num_results}
\input{sec7-conclusion}

\if1\blind{
\section*{Acknowledgments}
\noindent
Glaws' work is supported by the Ben L. Fryrear Ph.D. Fellowship in Computational Science at the Colorado School of Mines and the Department of Defense, Defense Advanced Research Project Agency's program Enabling Quantification of Uncertainty in Physical Systems. Constantine's work is supported by the US Department of Energy Office of Science, Office of Advanced Scientific Computing Research, Applied Mathematics program under Award Number DE-SC-0011077.
}\fi

\bibliographystyle{natbib}
\bibliography{lanczos-paper}

\end{document}

%% file: pauls-macros.tex
\graphicspath{{./}{./figs/}}

  \def\clap#1{\hbox to 0pt{\hss#1\hss}}

\providecommand{\mat}[1]{\bm{#1}}%
\renewcommand{\vec}[1]{\mathbf{#1}}
\newcommand{\vecalt}[1]{\bm{#1}}

\providecommand{\mA}{\ensuremath{\mat{A}}}

\providecommand{\mC}{\ensuremath{\mat{C}}}

\providecommand{\mE}{\ensuremath{\mat{E}}}

\providecommand{\mH}{\ensuremath{\mat{H}}}
\providecommand{\mI}{\ensuremath{\mat{I}}}
\providecommand{\mJ}{\ensuremath{\mat{J}}}

\providecommand{\mP}{\ensuremath{\mat{P}}}
\providecommand{\mQ}{\ensuremath{\mat{Q}}}

\providecommand{\mT}{\ensuremath{\mat{T}}}

\providecommand{\mV}{\ensuremath{\mat{V}}}
\providecommand{\mW}{\ensuremath{\mat{W}}}

\providecommand{\mLambda}{\ensuremath{\mat{\Lambda}}}
\providecommand{\mSigma}{\ensuremath{\mat{\Sigma}}}

\providecommand{\mzero}{\ensuremath{\mat{0}}}

\providecommand{\va}{\ensuremath{\vec{a}}}

\providecommand{\ve}{\ensuremath{\vec{e}}}

\providecommand{\vg}{\ensuremath{\vec{g}}}

\providecommand{\vu}{\ensuremath{\vec{u}}}
\providecommand{\vv}{\ensuremath{\vec{v}}}
\providecommand{\vw}{\ensuremath{\vec{w}}}
\providecommand{\vx}{\ensuremath{\vec{x}}}

\providecommand{\vmu}{\ensuremath{\vecalt{\mu}}}
\providecommand{\vphi}{\ensuremath{\vecalt{\phi}}}

\newcommand{\sF}{\mathcal{F}}
\newcommand{\sG}{\mathcal{G}}
\newcommand{\sI}{\mathcal{I}}

\newcommand{\sX}{\mathcal{X}}

\newcommand{\Exp}[1]{\mathbb{E}\left[#1\right]}

\newcommand{\Cov}[1]{\mathbb{C}\operatorname{ov}\left[#1\right]}

\newcommand{\diag}[1]{\operatorname{diag}\left(#1\right)}

\newcommand{\bmat}[1]{\begin{bmatrix}#1\end{bmatrix}}

\newcommand{\CIR}{\mC_{\operatorname{IR}}}
\newcommand{\hCIR}{\hat{\mC}_{\operatorname{IR}}}
\newcommand{\CSIR}{\mC_{\operatorname{SIR}}}
\newcommand{\hCSIR}{\hat{\mC}_{\operatorname{SIR}}}
\newcommand{\CAVE}{\mC_{\operatorname{AVE}}}
\newcommand{\hCAVE}{\hat{\mC}_{\operatorname{AVE}}}
\newcommand{\CSAVE}{\mC_{\operatorname{SAVE}}}
\newcommand{\hCSAVE}{\hat{\mC}_{\operatorname{SAVE}}}

%% file: sec1-intro.tex
\section{Introduction and background}
\label{sec:intro}

Increases in computational power have enabled the simulation of complex physical processes by models that more accurately reflect reality. Scientific studies that employ such models are often referred to as \emph{computer experiments}~\citep{Sacks89, Koehler96}. In this context, the computer simulation is represented by a deterministic function that maps the simulation inputs to some outputs of interest. These functions rarely have closed-form expressions and are often expensive to evaluate. To enable computer experiments with expensive models, we may construct a cheaper surrogate---e.g., a response surface approximation~\citep{Myers1995}---that can be evaluated quickly. However, building a surrogate suffers from \emph{the curse of dimensionality}~\citep{Traub1998,Donoho2000}; loosely speaking, the number of model evaluations required to achieve a prescribed accuracy grows exponentially with the input space dimension.

One approach to combat the curse of dimensionality is to reduce the dimension of the input space. Sufficient dimension reduction (SDR) provides a theoretical framework for dimension reduction in regression problems~\citep{Cook98, Ma13} and has been used as a method for dimension reduction in deterministic functions~\citep{Cook94b, Li16, Glaws17}. There are a variety of methods for SDR including ordinary least squares (OLS)~\citep{OLS89}, principal Hessian directions (pHd)~\citep{Li92}, among others~\citep{Cook05, Li07, Cook09}. In this paper, we consider two of the earliest methods introduced for SDR---the sliced inverse regression (SIR)~\citep{Li91} and the sliced average variance estimation (SAVE)~\citep{Cook91}.

SIR and SAVE each approximate a specific matrix of expectations using a slice-based approach. In the context of deterministic functions, these expectations become Lebesgue integrals, and the slicing can be interpreted as a Riemann sum approximation of the integrals~\citep{Davis84}. The approximation quality of these matrices depends on the number of terms in the Riemann sum (or slices) used. We introduce new algorithms---Lanczos-Stieltjes inverse regression (LSIR) and Lanczos-Stieltjes average variance estimation (LSAVE)---that improve accuracy and convergence rates compared to slicing by employing well-known tools from numerical analysis. The new algorithms enable use of high-order quadrature approximation when the number of inputs is sufficiently small to permit tensor product extensions to univariate quadrature---say, fewer than ten. Furthermore, the Lanczos-Stieltjes approaches perform as well as the best-case scenario of the slice-based algorithms when using Monte Carlo sampling as the basis for multivariate integration. This extends the methodology to regression problems and problems of very large dimension.

The remainder of this paper is structured as follows. In Section \ref{sec:IR_methods}, we review SDR theory and the SIR and SAVE algorithms. We examine the interpretation of these methods in the context of computer experiments in Section \ref{sec:comp_exp}. In Section \ref{sec:num_tools}, we review key elements of classical numerical analysis that we exploit in the new Lanczos-Stieltjes algorithms. We introduce and discuss the new algorithms in Section \ref{sec:LS_methods}. Finally, in Section \ref{sec:numericalresults}, we numerically study various aspects of the newly proposed algorithms on several test problems and compare the results to those from the traditional SIR and SAVE algorithms.

%% file: sec2-IR_methods.tex
\section{Inverse regression methods}
\label{sec:IR_methods}

Sufficient dimension reduction (SDR) enables dimension reduction in regression problems. In this section, we provide a brief overview of SDR and two methods for its implementation. More detailed discussions of SDR can be found in~\citet{Cook98} and~\citet{Ma13}.

The regression problem considers a given set of predictor/response pairs, denoted by
\begin{equation}
\label{eq:reg_prob}
\{ [ \, \vx_i^T \, , \, y_i \,, ] \} , \quad i = 0,\dots,N-1,
\end{equation}
where $\vx_i \in \mathbb{R}^m$ is the vector-valued predictor and $y_i \in \mathbb{R}$ is the scalar-valued response. These pairs are assumed to be drawn according to an unknown joint distribution $\pi_{\vx, y}$. The goal of regression is to approximate statistical properties (e.g., the CDF, expectation, variance) of the conditional random variable $y|\vx$ using the predictor/response pairs. However, such characterization of $y|\vx$ becomes difficult when the dimension of the predictor space $m$ becomes large. SDR aids in enabling the study of large-dimensional regression problems by seeking $\mA \in \mathbb{R}^{m \times n}$ with $n < m$ such that 
\begin{equation}
\label{eq:same_CDF}
y | \vx \,\sim\, y|\mA^T \vx ,
\end{equation}
provided such an $\mA$ exists. That is, $y|\vx$ and $y|\mA^T \vx$ have the same distribution so that no statistical information is lost by studying the response conditioned on the reduced predictors rather than the original predictors~\citep[Ch. 6]{Cook98}. In this sense, SDR can reduce the dimension of the regression problem.

We consider two methods for SDR---sliced inverse regression (SIR)~\citep{Li91} and sliced average variance estimation (SAVE)~\citep{Cook91}. These methods fall into a class of techniques known as \emph{inverse regression methods}, named as such because they search for dimension reduction in $y|\vx$ by studying the \emph{inverse regression} $\vx|y$~\citep{Adragni09}. The inverse regression is an $m$-dimensional conditional random variable parameterized by the scalar-valued response. This replaces the single $m$-dimensional problem with $m$ one-dimensional problems. The SIR and SAVE algorithms construct matrices using statistical characteristics of $\vx|y$ whose column spaces approximate $\mA$'s from \eqref{eq:same_CDF}.

Before describing the SIR and SAVE algorithms, we introduce the assumption of standardized predictors,
\begin{equation}
\label{eq:stand_preds}
\Exp{\vx} \;=\; \mzero
\quad \text{and} \quad
\Cov{\vx} \;=\; \mI .
\end{equation}
This assumption simplifies discussion and does not restrict the SDR problem~\citep[Ch. 11]{Cook98}. We assume \eqref{eq:stand_preds} holds for the remainder of this section.

Sliced inverse regression is based on the covariance of the expectation of the inverse regression,
\begin{equation}
\label{eq:CIR_reg}
\CIR \;=\; \Cov{\Exp{\vx|y}} .
\end{equation}
We denote this matrix by $\CIR$ to emphasize that it does not include the ``sliced'' component of the SIR algorithm. Sliced average variance estimation is based on the matrix
\begin{equation}
\label{eq:CAVE_reg}
\CAVE \;=\;\Exp{\left( \mI - \Cov{\vx|y} \right)^2} ,
\end{equation}
which again does not include the sliced component of the SAVE algorithm. The slicing method in SIR and SAVE enables numerical approximation of \eqref{eq:CIR_reg} and \eqref{eq:CAVE_reg}. Let $J_r$ for $r = 0,\dots,R-1$ denote $R$ intervals that partition the range of response values into slices. Note that this partitioning can include semi-infinite intervals if the range of response values is unbounded. Define the sliced mapping of the output
\begin{equation}
\label{eq:hy_eq_r}
h(y) \;=\; r
\quad \text{where} \quad
y \;\in\; J_r .
\end{equation}
The SIR and SAVE algorithms consider dimension reduction relative to the conditional random variable $h(y)|\vx$, not $y|\vx$. That is, these methods approximate the matrices
\begin{equation}
\begin{aligned}
\CSIR \;&=\; \Cov{\Exp{\vx|h(y)}}, \\
\CSAVE \;&=\; \Exp{\left( \mI - \Cov{\vx|h(y)} \right)^2} ,
\end{aligned}
\end{equation}
respectively, where the notation $\CSIR$ and $\CSAVE$ indicate that they are in terms of the sliced function $h(y)$. Under mild assumptions on the predictor space, the column spaces of these matrices have been shown to estimate the column space of $\mA$ from \eqref{eq:same_CDF}~\citep[Ch. 11]{Cook98}~\citep{Cook00}.

Algorithm \ref{alg:SIR} provides an outline for the SIR algorithm while SAVE is contained in Algorithm \ref{alg:SAVE}. Note that these algorithms typically end with an eigendecomposition of the sample estimates $\hCSIR$ and $\hCSAVE$ of the population statistics $\CSIR$ and $\CSAVE$. The eigenvectors of the sample estimates approximate $\mA$ from \eqref{eq:same_CDF}. However, for the purposes of this paper, our interest in these algorithms is how we can interpret them as approximations of $\CIR$ and $\CAVE$ from \eqref{eq:CIR_reg} and \eqref{eq:CAVE_reg}.

\begin{algorithm}
\caption{Sliced inverse regression (SIR)} \label{alg:SIR}
\textbf{Given:} predictor/response pairs $\{ [ \, \vx_i^T \, , \, y_i \, ] \}$, $i = 0,\dots,N-1$ drawn according to $\pi_{\vx, y}$ \\
\textbf{Assumptions:} The predictors are standardized such that \eqref{eq:stand_preds} is satisfied. \vspace{-0.8em}
\begin{enumerate}
\item Define a sliced partition of the observed response space, $J_r$ for $r = 0,\dots,R-1$. Let $\sI_r \subset \{ 0,\dots,N-1 \}$ be indices $i$ for which $y_i \in J_r$ and $N_r$ be the cardinality of $\sI_r$.
\item For $r = 0,\dots,R-1$, compute the in-slice sample expectation
\begin{equation}
\hat{\vmu} (r) \;=\; \frac{1}{N_r} \sum_{i \in \sI_r} \vx_i .
\end{equation} \vspace*{-2.5em}
\item Compute the sample matrix
\begin{equation}
\label{eq:hDSIR}
\hCSIR \;=\; \frac{1}{N} \sum_{r=0}^{R-1} N_r \, \hat{\vmu} (r) \, \hat{\vmu} (r)^T .
\end{equation} \vspace*{-1em}
\end{enumerate}
\end{algorithm}

\begin{algorithm}
\caption{Sliced average variance estimation (SAVE)} \label{alg:SAVE}
\textbf{Given:} predictor/response pairs $\{ [ \, \vx_i^T \, , \, y_i \, ] \}$, $i = 0,\dots,N-1$ drawn according to $\pi_{\vx, y}$ \\
\textbf{Assumptions:} The predictors are standardized such that \eqref{eq:stand_preds} is satisfied. \vspace{-0.8em}
\begin{enumerate}
\item Define a sliced partition of the observed response space, $J_r$ for $r = 0,\dots,R-1$. Let $\sI_r \subset \{ 0,\dots,N-1 \}$ be indices $i$ for which $y_i \in J_r$ and $N_r$ be the cardinality of $\sI_r$.
\item For $r = 0,\dots,R-1$, \vspace{-0.8em}
\begin{enumerate}
\item Compute the in-slice sample expectation
\begin{equation}
\hat{\vmu} (r) \;=\; \frac{1}{N_r} \sum_{i \in \sI_r} \vx_i .
\end{equation} \vspace*{-2.5em}
\item Compute the in-slice sample covariance
\begin{equation}
\hat{\mSigma} (r) \;=\; \frac{1}{N_r - 1} \sum_{i \in \sI_r} \left( \vx_i - \hat{\vmu} (r) \right) \left( \vx_i - \hat{\vmu} (r) \right)^T .
\end{equation} \vspace*{-2.5em}
\end{enumerate}
\item Compute the sample matrix
\begin{equation}
\label{eq:hDSAVE}
\hCSAVE \;=\; \frac{1}{N} \sum_{r=0}^{R-1} N_r \left( \mI - \hat{\mSigma} (r) \right)^2 .
\end{equation} \vspace*{-2em}
\end{enumerate}
\end{algorithm}

The sliced mapping of the response enables computation of the sample estimates
\begin{equation} \label{eq:sample_est}
\hat{\vmu} (r) \;\approx\; \Exp{\vx|h(y)} 
\quad \text{and} \quad
\hat{\mSigma} (r) \;\approx\; \Cov{\vx|h(y)} 
\end{equation}
by binning the response data. The approximation in \eqref{eq:sample_est} depends on the number of predictor/response pairs $N$. The sample estimates of $\CSIR$ and $\CSAVE$ have been shown to be $N^{-1/2}$ consistent~\citep{Li91, Cook00}. Accuracy of $\CSIR$ and $\CSAVE$ in approximating $\CIR$ and $\CAVE$, respectively, depends on the slicing applied to the response space. Defining the slices $J_r$, $r = 0,\dots,R-1$ such that they each contain approximately equal numbers of samples can improve the approximation by balancing accuracy in the sample estimates \eqref{eq:sample_est} across all slices~\citep{Li91}. By increasing the number $R$ of slices, we increase accuracy due to slicing but may reduce accuracy of \eqref{eq:sample_est} by decreasing the number of samples within each slice~\citep[Ch. 6]{Cook98}. We discuss the idea of convergence in terms of the number of slices in more detail in Section \ref{sec:comp_exp}. Generally, the slicing approach in Algorithms \ref{alg:SIR}and \ref{alg:SAVE} is justified by properties of SDR that relate dimension reduction of $h(y)|\vx$ to dimension reduction of $y|\vx$~\citep[Ch. 6]{Cook98}. 

%% file: sec3-comp_exp.tex
\section{SIR \& SAVE for deterministic functions}
\label{sec:comp_exp}

In this section, we consider the application of SIR and SAVE to deterministic functions, such as those arising in computer experiments. Computer experiments employ models of real-world phenomena, which we represent as deterministic functions that map $m$ simulation inputs to a scalar-valued output~\citep{Sacks89, Koehler96}. We write this as
\begin{equation}
\label{eq:y_eq_fx}
y \;=\; f(\vx) ,
\qquad \vx \in \sX \subseteq \mathbb{R}^m , \quad y \in \sF \subseteq \mathbb{R},
\end{equation}
where $f: \sX \rightarrow \sF$ is assumed measurable and the random variables $\vx$, $y$ are the function inputs and output, respectively. Let $\rho$ be the probability measure induced by $\vx$ over $\sX$. Unlike the joint distribution $\pi_{\vx, y}$, this measure is assumed to be known. Furthermore, we may choose the values of $\vx$ at which we wish to evaluate \eqref{eq:y_eq_fx} according to $\rho$. We define the probability measure induced by $y$ over $\sF$ to be $\gamma$. This measure is fully determined by $\rho$ and $f$, but its form is assumed to be unknown. We can estimate $\gamma$ by point evaluations of $f$.

\begin{remark}
In the context of deterministic functions, SDR may be viewed as a method for ridge recovery~\citep{Glaws17}. That is, it searches for $\mA \in \mathbb{R}^{m \times n}$ with $n < m$ assuming that
\begin{equation}
\label{eq:y_eq_fx_eq_gAx}
y \;=\; f(\vx)
\;=\; g(\mA^T \vx)
\end{equation}
for some $g : \mathbb{R}^n \rightarrow \mathbb{R}$. When \eqref{eq:y_eq_fx_eq_gAx} holds for some $g$ and $\mA$, $f(\vx)$ is called a \emph{ridge function}~\citep{Pinkus15}. \citet{Cook94b} and~\citet{Li16} apply SDR for ridge recovery in deterministic functions.
\end{remark}
 
In Section \ref{sec:LS_methods}, we introduce new algorithms for estimating the $\CIR$ and $\CAVE$ matrices which form the basis of the SIR and SAVE algorithms. First, we must provide context for applying these methods to deterministic functions. We start by defining an assumption over the input space $\sX$.
\begin{assume}
\label{assume:stand_inputs}
The inputs have finite fourth moments and are standardized such that
\begin{equation}
\label{eq:stand_inputs}
\int \vx \, d \rho (\vx) \;=\; \mzero
\quad \text{and} \quad
\int \vx \, \vx^T \, d \rho (\vx) \;=\; \mI .
\end{equation}
\end{assume}
Assumption \ref{assume:stand_inputs} is similar to the assumption of standardized predictors in \eqref{eq:stand_preds}. However, it also assumes finite fourth moments of the predictors. This assumption is important for the new algorithms introduced in Section \ref{sec:LS_methods}. We assume Assumption \ref{assume:stand_inputs} holds throughout this paper.

Recall that SIR and SAVE employ the inverse regression $\vx|y$. For deterministic functions, this is the inverse image of the function $f$,
\begin{equation}
f^{-1} (y) \;=\; \left\{ \vx \;\in\; \sX \, : \, f(\vx) \;=\; y \right\} ,
\end{equation}
and we define the conditional measure $\sigma_{\vx|y}$ as the restriction of $\rho$ to $f^{-1}$~\citep{Chang97}. We denote this measure by $\sigma_{\vx|y}$ to emphasize its connection to the inverse regression $\vx|y$. 


The matrices $\CIR$ and $\CAVE$ contain the conditional expectation and the conditional covariance of the inverse regression (see \eqref{eq:CIR_reg} and \eqref{eq:CAVE_reg}, respectively). These statistical quantities can be expressed as integrals with respect to the conditional measure $\sigma_{\vx|y}$,
\begin{equation}
\label{eq:int_forms1}
\begin{aligned}
\vmu (y) \;&=\; \int \vx \, d \sigma_{\vx|y} (\vx) , \\
\mSigma (y) \;&=\; \int \left( \vx - \vmu (y) \right) \, \left( \vx - \vmu (y) \right)^T \, d \sigma_{\vx|y} (\vx) .
\end{aligned}
\end{equation}
We denote the conditional expectation and covariance by $\vmu (y)$ and $\mSigma (y)$, respectively, to emphasize that these quantities are functions of the output. Using \eqref{eq:int_forms1}, we can express the $\CIR$ and $\CAVE$ matrices as integrals over the output space $\sF$,
\begin{equation}
\label{eq:CIR_CAVE}
\begin{aligned}
\CIR \;&=\; \int \vmu (y) \, \vmu (y)^T \, d \gamma (y) , \\
\CAVE \;&=\: \int \left( \mI - \mSigma (y) \right)^2 \, d \gamma (y) .
\end{aligned}
\end{equation}

The integrals in \eqref{eq:int_forms1} are difficult to approximate due to the complex structure of $\sigma_{\vx|y}$. Algorithms \ref{alg:SIR} and \ref{alg:SAVE} slice the range of output values to enable this approximation. Let $h(y)$ denote the slicing function from \eqref{eq:hy_eq_r} and define the probability mass function
\begin{equation}
\omega (r)
\;=\;
\int_{J_r} \, d \gamma (y) ,
\qquad r \;=\; 0,\dots,R-1 ,
\end{equation}
where $J_r$ is the $r$th interval over the range of outputs. We can then express the slice-based matrices $\CSIR$ and $\CSAVE$ as sums over the $R$ slices,
\begin{equation}
\label{eq:CSIR_CSAVE}
\begin{aligned}
\CSIR \;&=\; \sum_{r=0}^{R-1} \omega (r) \, \vmu (r) \, \vmu (r)^T , \\
\CSAVE \;&=\: \sum_{r=0}^{R-1} \omega (r) \, \left( \mI - \mSigma (r) \right)^2 ,
\end{aligned}
\end{equation}
where $\vmu (r)$ and $\mSigma (r)$ are as in \eqref{eq:int_forms1} but applied to the sliced output $h(y) = r$.

This exercise of expressing the various elements of SIR and SAVE as integrals emphasizes the two levels of approximation occurring in these algorithms:  (i) approximation in terms of the number of samples $N$ and (ii) approximation in terms of the number of slices $R$. That is,
\begin{equation}
\label{eq:2lvl_approx}
\begin{gathered}
\hCSIR \;\approx\; \CSIR \;\approx\; \CIR , \\
\hCSAVE \;\approx\; \CSAVE \;\approx\; \CAVE .
\end{gathered}
\end{equation}
As mentioned in Section \ref{sec:IR_methods}, approximation due to sampling (i.e., the leftmost approximation in \eqref{eq:2lvl_approx}) has been shown to be $N^{-1/2}$ consistent~\citep{Li91, Cook00}. From the integral perspective, the slicing approach can be interpreted as a Riemann sum approximation of the integrals in \eqref{eq:CIR_CAVE}. Riemann sum approximations estimate integrals by a finite summation of values. These approximations converge like $R^{-1}$ for continuous functions, where $R$ denotes the number of Riemann sums (or slices)~\citep[Ch. 2]{Davis84}. 

In the next section, we introduce several numerical tools and link them to the various elements of the SIR and SAVE algorithms discussed so far. Ultimately, we use these tools, including orthogonal polynomials and numerical quadrature, in the proposed algorithms to enable approximation of $\CIR$ and $\CAVE$ without using Riemann sums. 

%% file: sec4-num_tools.tex
\section{Orthogonal polynomials and Gauss-Christoffel quadrature}
\label{sec:num_tools}

In this section, we review elements of orthogonal polynomials and numerical quadrature. These topics are fundamental to the field of numerical analysis and have been studied extensively; detailed discussions are available in~\citet{Liesen13, Gautschi04, Meurant06, Golub10}. Our discussion is based on these references; however, we limit it to the key concepts necessary to develop the new algorithms for approximating the matrices in \eqref{eq:CIR_CAVE}. We begin with the Stieltjes procedure for constructing orthonormal polynomials with respect to a given measure. We then relate this procedure to Gauss-Christoffel quadrature, polynomial expansions, and the Lanczos algorithm.


\subsection{The Stieltjes procedure}
\label{subsec:stieltjes}

The Stieltjes procedure recursively constructs a sequence of polynomials that are orthonormal with respect to a given measure~\citep[Ch. 3]{Liesen13}. Let $\gamma$ denote a given probability measure over $\mathbb{R}$, and let $\phi , \psi : \mathbb{R} \rightarrow \mathbb{R}$ be two scalar-valued functions that are square-integrable with respect to $\gamma$. The continuous inner product relative to $\gamma$ is
\begin{equation} \label{eq:inner_prod}
( \phi , \psi )
\;=\;
\int \phi(y) \, \psi(y) \, d \gamma (y) ,
\end{equation}
and the induced norm is $|| \phi || = \sqrt{( \phi , \phi )}$. A sequence of polynomials $\{ \phi_0 , \phi_1 , \phi_2 , \dots \}$ is orthonormal with respect to $\gamma$ if
\begin{equation}
( \phi_i , \phi_j )
\;=\;
\delta_{i,j} ,
\quad i,j = 0,1,2,\dots ,
\end{equation}
where $\delta_{i,j}$ is the Kronecker delta. Algorithm \ref{alg:Stieltjes} contains a method for constructing a sequence of orthonormal polynomials $\{ \phi_0 , \phi_1 , \phi_2 , \dots \}$ relative to $\gamma$. This algorithm is known as the Stieltjes procedure and was first introduced in~\citet{Stieltjes84}. 

\begin{algorithm}
\caption{Stieltjes procedure~\citep[Section 2.2.3.1]{Gautschi04}} \label{alg:Stieltjes}
\textbf{Given:} probability measure $\gamma$ \\
\textbf{Assumptions:} Let $\phi_{-1} (y) = 0$ and $\tilde{\phi}_1 (y) = 1$. \\
\textbf{for } $i \;=\; 0, 1, 2, \dots$ \\
\hspace*{1em} $\beta_i \;=\; ||\tilde{\phi}_i||$ \\
\hspace*{1em} $\phi_i \;=\; \tilde{\phi}_i \, / \, \beta_i$ \\
\hspace*{1em} $\alpha_i \;=\; (y \, \phi_i, \phi_i)$ \\
\hspace*{1em} $\tilde{\phi}_{i+1} \;=\; (y - \alpha_i) \phi_i - \beta_i \phi_{i-1}$ \\
\textbf{end}
\end{algorithm}

In the last step of Algorithm \ref{alg:Stieltjes}, we see the three-term recurrence relationship for orthonormal polynomials,
\begin{equation} \label{eq:three_term_rec}
\beta_{i+1} \phi_{i+1} (y)
\;=\; 
(y - \alpha_i) \phi_i (y) - \beta_i \phi_{i-1} (y) ,
\qquad i = 0, 1, 2, \dots .
\end{equation}
Any sequence of polynomials that satisfies \eqref{eq:three_term_rec} is orthonormal with respect to the given measure. If we consider the first $k$ terms, then we can rearrange it to obtain
\begin{equation} \label{eq:three_term_rec2}
y \, \phi_i (y)
\;=\;
\beta_i \phi_{i-1} (y) + \alpha_i \phi_i (y) + \beta_{i+1} \phi_{i+1} (y) ,
\qquad i = 0, 1, 2, \dots, k-1 .
\end{equation}
Let $\vphi (y) = [ \, \phi_0 (y) \, , \, \phi_1 (y) \, , \, \dots \, , \, \phi_{k-1} (y) \, ]^T$. We can then write \eqref{eq:three_term_rec2} in matrix form as
\begin{equation} \label{eq:3_term_matrix}
y \, \vphi (y)
\;=\;
\mJ \, \vphi (y) + \beta_k \, \phi_k (y) \, \ve_k ,
\end{equation}
where $\ve_k \in \mathbb{R}^k$ is the vector of zeros with a one in the last entry and $\mJ \in \mathbb{R}^{k \times k}$ is
\begin{equation} \label{eq:J_Jacobi}
\mJ
\;=\;
\bmat{
\alpha_0 & \beta_1 & & & \\
\beta_1 & \alpha_1 & \beta_2 & & \\
 & \ddots & \ddots & \ddots & \\
 & & \beta_{k-2} & \alpha_{k-2} & \beta_{k-1} \\
 & & & \beta_{k-1} & \alpha_{k-1}} ,
\end{equation}
where $\alpha_i$, $\beta_i$ are the recurrence coefficients from Algorithm \ref{alg:Stieltjes}. Note that this matrix---known as the Jacobi matrix---is symmetric and tridiagonal. Let the eigendecomposition of $\mJ$ be
\begin{equation}
\mJ
\;=\;
\mQ \mLambda \mQ^T .
\end{equation}
From \eqref{eq:3_term_matrix}, the eigenvalues of $\mJ$, denoted by $\lambda_i$, $i = 0,\dots,k-1$, are the zeros of the degree-$k$ polynomial $\phi_k (y)$. Furthermore, the eigenvector associated with $\lambda_i$ is $\vphi (\lambda_i)$. We assume the eigenvectors of $\mJ$ are normalized such that $\mQ$ is an orthogonal matrix with entries
\begin{equation} \label{eq:Q_ij}
\left( \mQ \right)_{i,j}
\;=\;
\frac{\phi_i (\lambda_j)}{|| \vphi (\lambda_j) ||_2} ,
\qquad i,j = 0,\dots,k-1 ,
\end{equation}
where the notation $(\cdot)_{i,j}$ denotes the element in the $i$th row and $j$th column of the given matrix. 

We end this section with brief a note about Fourier expansion of functions in terms of orthonormal polynomials. If a given function $g(y)$ is square integrable with respect to $\gamma$, then $g$ admits a mean-squared convergent Fourier series in terms of the orthonormal polynomials,
\begin{equation} \label{eq:fourier_exp}
g(y)
\;=\;
\sum_{i=0}^{\infty} g_i \, \phi_i (y) ,
\end{equation}
for $y$ in the support of $\gamma$ and where equality is in the $L_2 (\gamma)$ sense. By orthogonality of the polynomials, the Fourier coefficients $g_i$ are
\begin{equation} \label{eq:fourier_coeffs}
g_i
\;=\;
\left( g , \phi_i \right) .
\end{equation}
This polynomial approximation plays an important role in the algorithms introduced in Section \ref{sec:LS_methods}.


\subsection{Gauss-Christoffel quadrature}
\label{subsec:gauss-christoffel}

We next discuss the Gauss-Christoffel quadrature for numerical integration and show its connection to orthonormal polynomials~\citep[Ch. 3]{Liesen13}. Given a measure $\gamma$ and integrable function $g(y)$, a $k$-point quadrature rule approximates the integral of $g$ with respect to $\gamma$ by a weighted sum of $g$ evaluated at $k$ input values,
\begin{equation}
\int g(y) \, d \gamma (y)
\;=\;
\sum_{i = 0}^{k-1} \omega_i \, g(\lambda_i) + r_k .
\end{equation}
The $\lambda_i$'s are the quadrature nodes and the $\omega_i$'s are the associated quadrature weights. The $k$-point quadrature approximation error is contained in the residual term $r_k$. We can minimize $|r_k|$ by choosing the quadrature nodes and weights appropriately. The nodes and weights of the Gauss-Christoffel quadrature maximize the polynomial \emph{degree of exactness}, which means the quadrature rules exactly integrate (i.e., $r_k = 0$) all polynomials of degree at most $k$. The $k$-point Gauss-Christoffel quadrature rule has polynomial degree of exactness $2k - 1$. Furthermore, the Gauss-Christoffel quadrature has been shown to converge exponentially at a rate $\rho^{-k}$ when the integrand is analytic~\citep[Ch. 19]{Trefethen13}. The base $\rho > 1$ relates to the size of the function's domain of analytical continuability For functions with $p-1$ continuous derivatives, the Gauss-Christoffel quadrature converges like $k^{-(2p + 1)}$.

The Gauss-Christoffel quadrature nodes and weights depend on the given measure $\gamma$. They can be obtained through the eigendecomposition of $\mJ$~\citep{Golub69}. Recall from \eqref{eq:J_Jacobi} that $\mJ$ is the matrix of recurrence coefficients resulting from $k$ steps of the Stieltjes procedure. The eigenvalues of $\mJ$ are the zeros of the $k$-degree orthonormal polynomial $\phi_k (y)$. These zeros are the nodes of the $k$-point Gauss-Christoffel quadrature rule with respect to $\gamma$ (i.e., $\phi_k (\lambda_i) = 0$ for $i = 0,\dots,k-1$). The associated weights are the squares of the first entry of each normalized eigenvector,
\begin{equation}
\omega_i
\;=\;
\left( \mQ \right)_{0,i}^2
\;=\;
\frac{1}{||\vphi (\lambda_i)||_2^2} ,
\quad i = 0, \dots ,k-1.
\end{equation}

Recall that the Stieltjes procedure employs the inner product from \eqref{eq:inner_prod} to define the recurrence coefficients $\alpha_i$, $\beta_i$. In the next section, we consider the Lanczos algorithm and explore conditions under which it be may considered a discrete analog to the Stieltjes procedure in Section \ref{subsec:stieltjes}. Before we make this connection, we define the discrete inner product as the Gauss-Christoffel quadrature approximation of the continuous inner product, 
\begin{equation}
\label{eq:disc_inner_prod}
\left( \phi , \psi \right)_k
\;=\;
\sum_{i = 0}^{k-1} \omega_i \, \phi (\lambda_i) \, \psi (\lambda_i)
\;\approx\;
\left( \phi , \psi \right) ,
\end{equation}
where the subscript $k$ denotes the order of the quadrature rule used. Similarly, we define the discrete norm as $|| \phi ||_k = \sqrt{( \phi , \phi )_k}$. Note that we could use any numerical integration technique to define a discrete inner product. This fact appears in the numerical experiments in Section \ref{sec:numericalresults}; however, for the discussion in this section and in Section \ref{sec:LS_methods} we define the discrete inner product as in \eqref{eq:disc_inner_prod} using the Gauss-Christoffel quadrature.

The pseudospectral expansion approximates the Fourier expansion from \eqref{eq:fourier_exp} by truncating the series after $k$ terms and approximating the Fourier coefficients in \eqref{eq:fourier_coeffs} using the discrete inner product~\citep{Constantine12b}. We write this series for a given square-integrable function $g(y)$ as
\begin{equation} \label{eq:pseudo_exp}
\hat{g} (y)
\;=\;
\sum_{i=0}^{k-1} \hat{g}_i \, \phi_i (y) ,
\end{equation}
where the pseudospectral coefficients are
\begin{equation} \label{eq:pseudo_coeffs}
\hat{g}_i
\;=\;
\left( g , \phi_i \right)_k .
\end{equation}
Note that the approximation of $g(y)$ by $\hat{g}(y)$ depends on two levels of approximation: (i) the accuracy of the approximated coefficients $\hat{g}_i$, $i = 0,\dots,k-1$ and (ii) the magnitude of the coefficients in the omitted terms $g_i$, $i = k,k+1,\dots$. Using a higher-order quadrature rule improves (i) and including more terms in the truncated series improves (ii). The pseudospectral approximation and this two-level convergence play an important role in the new algorithms proposed in Section \ref{sec:LS_methods}.


\subsection{The Lanczos algorithm}
\label{subsec:lanczos}

The Lanczos algorithm was originally introduced as an iterative scheme for approximating eigenvalues and eigenvectors of linear differential operators~\citep{Lanczos50}. Given a symmetric $N \times N$ matrix $\mA$, it constructs a symmetric, tridiagonal $k \times k$ matrix $\mT$ whose eigenvalues approximate those of $\mA$. Additionally, it produces an $N \times k$ matrix, denoted by $\mV = [ \, \vv_0 \, , \, \vv_1 \, , \dots , \, \vv_{k-1} \, ]$ where the $\vv_i$'s are the Lanczos vectors. These vectors transform the eigenvectors of $\mT$ into approximate eigenvectors of $\mA$. Algorithm \ref{alg:Lanczos} contains the steps of the Lanczos algorithm. Note that the inner products and norms in Algorithm \ref{alg:Lanczos} are the given by
\begin{equation}
\left( \vw , \vu \right) \;=\; \vw^T \vu ,
\qquad
|| \vw || \;=\; \sqrt{\left(\vw , \vw \right)} ,
\end{equation}
for vectors $\vw , \vu \in \mathbb{R}^N$. These do not directly relate to the discrete inner products introduced in Section \ref{subsec:gauss-christoffel}, although we make several connections later in this section.

\begin{algorithm}
\caption{Lanczos algorithm~\citep[Section 3.1.7.1]{Gautschi04}} \label{alg:Lanczos}
\textbf{Given:} A $N \times N$ symmetric matrix $\mA$. \\
\textbf{Assumptions:} Let $\vv_{-1} = \mzero \in \mathbb{R}^N$ and $\tilde{\vv}_0$ be an arbitrary nonzero vector of length $N$. \\
\textbf{for } $i \;=\; 0$ \textbf{ to } $k-1$, \\
\hspace*{1em} $\beta_i \;=\; || \tilde{\vv}_i ||$ \\
\hspace*{1em} $\vv_i \;=\; \tilde{\vv}_i \, / \, \beta_i$ \\
\hspace*{1em} $\alpha_i \;=\; \left( \mA \vv_i , \vv_i \right)$ \\
\hspace*{1em} $\tilde{\vv}_{i+1} \;=\; (\mA - \alpha_i \mI) \vv_i - \beta_{i-1} \vv_{i-1}$ \\
\textbf{end}
\end{algorithm}

After $k$ iterations, the Lanczos algorithm results in the recurrence relationship
\begin{equation} \label{eq:k_steps_Lanczos}
\mA \mV
\;=\;
\mV \mT + \beta_k \vv_k \ve_k^T ,
\end{equation}
where $\ve_k \in \mathbb{R}^k$ is the vector of zeros with a one in the last entry, $\mV$ is the matrix of Lanczos vectors, and $\mT$ is a symmetric, tridiagonal matrix of the recurrence coefficients, 
\begin{equation} \label{eq:T_Jacobi}
\mT
\;=\;
\bmat{
\alpha_0 & \beta_1 & & & \\
\beta_1 & \alpha_1 & \beta_2 & & \\
 & \ddots & \ddots & \ddots & \\
 & & \beta_{k-2} & \alpha_{k-2} & \beta_{k-1} \\
 & & & \beta_{k-1} & \alpha_{k-1}} .
\end{equation}
The matrix $\mT$ in \eqref{eq:T_Jacobi}, similar to $\mJ$ in \eqref{eq:J_Jacobi}, is the Jacobi matrix~\citep{Gautschi04}. The relationship between the matrices $\mJ$ and $\mT$ has been studied extensively~\citep{Gautschi02, Forsythe57, deBoor78}. 

We are interested in the use of the Lanczos algorithm as a discrete approximation to the Stieltjes procedure. Recall that Algorithm \ref{alg:Stieltjes} (the Stieltjes procedure) assumes a given measure $\gamma$, and consider the $N$-point Gauss-Christoffel quadrature rule with respect to $\gamma$, which has nodes and weights $\lambda_i$, $\omega_i$ for $i = 0,\dots,N-1$. This quadrature rule defines a discrete approximation of $\gamma$, which we denote by $\gamma_N$. Inner products with respect to $\gamma_N$ take the form of the discrete inner product in \eqref{eq:disc_inner_prod}. If we perform the Lanczos algorithm (Algorithm \ref{alg:Lanczos}) with inputs
\begin{equation}
\mA
\;=\;
\bmat{\lambda_0 & & \\ & \ddots & \\ & & \lambda_{N-1}} ,
\qquad
\tilde{\vv}_0
\;=\;
\bmat{\sqrt{\omega_0} \\ \vdots \\ \sqrt{\omega_{N-1}}} ,
\end{equation}
then the result is equivalent to running the Stieltjes procedure using this discrete inner product with respect to $\gamma_N$. Increasing the number of quadrature nodes $N$ improves the approximation of $\gamma_N$ to $\gamma$. It can be shown that the recurrence coefficients in the resulting Jacobi matrix will converge to the recurrence coefficients related to the Stieltjes procedure with respect to $\gamma$ as $N$ increases~\citep[Section 2.2]{Gautschi04}. In the next, section we show how this relationship between Stieltjes and Lanczos can be used to approximate composite functions, which is essential to understand the underpinnings of LSIR and LSAVE in Section \ref{sec:LS_methods}.


\subsection{Composite function approximation}
\label{subsec:comp_func}

The connection between Algorithms \ref{alg:Stieltjes} and \ref{alg:Lanczos} can be exploited for the approximation of composite functions~\citep{Constantine12a}. Consider a function of the form
\begin{equation}
\label{eq:hx_equals_gfx}
h(x)
\;=\;
g(f(x)) ,
\qquad x \in \sX \subseteq \mathbb{R}
\end{equation}
where
\begin{equation}
\label{eq:sX_sF_sG}
\begin{aligned}
f&: \sX \rightarrow \sF \subseteq \mathbb{R} , \\
g&: \sF \rightarrow \sG \subseteq \mathbb{R} .
\end{aligned}
\end{equation}
Assume the input space $\sX$ is weighted with a given probability measure $\rho$. This measure and $f$ induce a measure on $\sF$ that we denote by $\gamma$.  Note that the methodology described in this section can be extended to multivariate inputs (i.e., $\sX \subseteq \mathbb{R}^m$) through tensor product constructions and to multivariate outputs (i.e., $\sG \subseteq \mathbb{R}^n$) by considering each output individually. We need both of these extensions in Section \ref{sec:LS_methods}; however, we consider the scalar case here for clarity.

The goal here is to construct a pseudospectral expansion of $g$ using orthonormal polynomials with respect to $\gamma$ and a Gauss-Christoffel quadrature rule defined over $\sF$. In Section \ref{subsec:stieltjes}, we examined the Stieltjes procedure which constructs a sequence of orthonormal polynomials with respect to a given measure. In Section \ref{subsec:gauss-christoffel}, we showed this algorithm also produces the nodes and weights of the Gauss-Christoffel quadrature rule with respect to the given measure. In Section \ref{subsec:lanczos}, we saw how the Lanczos algorithm can be used to produce similar results for a discrete approximation of the given measure. All of this suggests a methodology for constructing a pseudospectral approximation of $g$. However, we constructed the discrete approximation to $\gamma$ in Section \ref{subsec:lanczos} using a quadrature rule relative to this measure. We cannot do this here since the form of $\gamma$ is unknown. We can construct the $N$-point Gauss-Christoffel quadrature rule on $\sX$ since $\rho$ is assumed known. Let $x_i$, $\nu_i$, $i = 0,\dots,N-1$ denote these quadrature nodes and weights. This quadrature rule defines a discrete approximation of $\rho$, which we write as $\rho_N$. We then approximate $\gamma$ by the discrete measure $\gamma_N$ by evaluating $f_i = f(x_i)$ for $i = 0,\dots,N-1$ and continue as before. Algorithm \ref{alg:Lanc_comp_funcs} outlines this process for obtaining orthonormal polynomials and a Gauss-Christoffel quadrature rule relative to $\gamma_N$.

\begin{algorithm}[!ht]
\caption{\citet{Constantine12a}} \label{alg:Lanc_comp_funcs}
\textbf{Given:} function $f$ and input probability measure $\rho$ \vspace*{-1em}
\begin{enumerate}
\item Obtain the $N$ nodes and weights for the Gauss-Christoffel quadrature rule with respect to $\rho$,
\begin{equation}
x_i , \; \nu_i ,
\quad \text{for} \quad i = 0,\dots,N-1 .
\end{equation}
\item Evaluate $f_i = f(x_i)$ for $i = 0, \dots ,N-1$.
\item Define
\begin{equation}
\label{eq:Lanc_for_comp}
\mA
\;=\;
\bmat{f_0 & & \\ & \ddots & \\ & & f_{N-1}} ,
\qquad
\tilde{\vv}_0
\;=\;
\bmat{\sqrt{\nu_0} \\ \sqrt{\nu_1} \\ \vdots \\ \sqrt{\nu_{N-1}}} .
\end{equation}
\item Perform $k$ iterations of the Lanczos algorithm (see Algorithm \ref{alg:Lanczos}) on $\mA$ with initial vector $\tilde{\vv}_0$ to obtain
\begin{equation}
\mA \mV
\;=\;
\mV \mT + \beta_k \vv_k \ve_k^T .
\end{equation}
\end{enumerate}
\end{algorithm}

Let the eigendecomposition of the resulting Jacobi matrix be
\begin{equation}
\mT
\;=\;
\mQ \mLambda \mQ^T .
\end{equation}
We know the eigenvalues of $\mT$ define the $k$-point Gauss-Christoffel quadrature nodes relative to $\gamma_N$. We denote these quadrature nodes by $\lambda^N_i$, $i = 0,\dots,k-1$, where the superscript $N$ indicates that these nodes are relative to the discrete measure $\gamma_N$. From \eqref{eq:Q_ij}, we know the normalized eigenvectors of $\mT$ have the form
\begin{equation}
\left( \mQ \right)_{i}
\;=\;
\frac{\vphi^N (\lambda_i^N)}{||\vphi^N (\lambda_i^N)||_2} ,
\qquad i = 0,\dots,k-1 ,
\end{equation}
where $\vphi^N (y) = [ \, \phi^N_0 (y) \, , \, \phi^N_1 (y) \, , \, \dots \, , \, \phi^N_{k-1} (y) \, ]^T$ and $\phi^N_i (y)$ denotes the $i$th orthonormal polynomial relative to $\gamma_N$. The quadrature weight associated with $\lambda^N_i$ is given by the square of the first element of the $i$th normalized eigenvector,
\begin{equation}
\omega^N_i
\;\approx\;
\left( \mQ \right)_{0, i}^2 ,
\qquad i = 0, \dots ,k-1 .
\end{equation}
From Section \ref{subsec:lanczos} we have convergence of the quadrature nodes $\lambda^N_i$ and weights $\omega^N_i$ to $\lambda_i$ and $\omega_i$, respectively, as $N$ increases. In this sense, we consider these quantities to be approximations of the quadrature nodes and weights relative to $\gamma$,
\begin{equation}
\label{eq:approx_quad}
\lambda^N_i \;\approx\; \lambda_i
\quad \text{and} \quad
\omega^N_i \;\approx\; \omega_i .
\end{equation}
In Section \ref{sec:numericalresults}, we numerically study this approximation.

An alternative perspective on this $k$-point Gauss-Christoffel quadrature rule is that of a discrete measure $\gamma_{N,k}$ that approximates the measure $\gamma_N$. That is,
\begin{equation}
\label{eq:gamma_approx}
\gamma_{N,k} \;\approx\; \gamma_N \;\approx\; \gamma .
\end{equation}
By taking more Lanczos iterations, we improve the leftmost approximation, and for $k = N$, we have $\gamma_{N,N} = \gamma_N$. By increasing $N$ (the number of quadrature nodes on $\sX$), we improve the rightmost approximation. This may be viewed as convergence of the discrete Lanczos algorithm to the continuous Stieltjes procedure. This mirrors the two-level approximation from \eqref{eq:2lvl_approx}. The leftmost approximation in each case is defined by a numerical integration rule; however, the quality of those integration rules is vastly different (Gauss-Christoffel quadrature in \eqref{eq:gamma_approx} versus Riemann sums in \eqref{eq:2lvl_approx}).

\citet{Constantine12a} show that the Lanczos vectors resulting from Algorithm \ref{alg:Lanc_comp_funcs} also contain useful information. The Lanczos vectors are of the form
\begin{equation}
\label{eq:lanczos_vectors}
(\mV)_{i,j}
\;\approx\;
\sqrt{\nu_i} \, \phi_j (f_i) ,
\qquad i = 0, \dots ,N-1 ,
\quad j = 0, \dots ,k-1 ,
\end{equation}
where $f_i$ is the evaluation of $f(\vx)$ at the quadrature node $\vx_i$ relative to $\rho$, $\nu_i$ is the associated quadrature weight, and $\phi_j$ is the $j$th-degree orthonormal polynomial with respect to $\gamma$. The approximation in \eqref{eq:lanczos_vectors} is due to the approximation of $\gamma$ by $\gamma_N$ and is in the same vein as the approximation in \eqref{eq:approx_quad}. We also numerically study this approximation in Section \ref{sec:numericalresults}.

The approximation method in~\citet{Constantine12a} suggests evaluating $g$ at the $k \ll N$ quadrature nodes resulting from Algorithm \ref{alg:Lanc_comp_funcs} and using these evaluations to construct a pseudospectral approximation. This allows for accurate estimation of $g$ while placing a majority of the computational cost on evaluating $f$ instead of both $f$ and $g$ in \eqref{eq:hx_equals_gfx}. Such an approach is valuable when $g$ is difficult to compute relative to $f$. In the next section, we explain how this methodology can be used to construct approximations to $\CIR$ and $\CAVE$ from \eqref{eq:CIR_CAVE}.

%% file: sec5-LS_methods.tex
\section{Lanczos-Stieltjes methods for inverse regression}
\label{sec:LS_methods}

In this section, we use the tools discussed in Section \ref{sec:num_tools} to develop a new Lanczos-Stieltjes approach to inverse regression methods. This approach avoids approximating $\CIR$ and $\CAVE$ by Riemann sums (or slicing) as discussed in Section \ref{sec:comp_exp}. Instead, we use orthonormal polynomials and quadrature approximations to build more accurate estimates of these matrices. For reference, we reiterate the problem setup from \eqref{eq:y_eq_fx},
\begin{equation}
\label{eq:y_eq_fx2}
y \;=\; f(\vx) ,
\qquad \vx \in \sX \subseteq \mathbb{R}^m , \quad y \in \sF \subseteq \mathbb{R},
\end{equation}
with $\sX$ weighted by the known probability measure $\rho$ and $\sF$ weighted by the unknown probability measure $\gamma$.


\subsection{Lanczos-Stieltjes inverse regression (LSIR)}
\label{subsec:LSIR}

In Section \ref{sec:comp_exp}, we showed that the SIR algorithm approximates the matrix
\begin{equation} \label{eq:CIR2}
\CIR
\;=\;
\int \vmu (y) \, \vmu (y)^T \, d \gamma (y)
\end{equation}
using a sliced mapping of the output (i.e., a Riemann sum approximation). We wish to approximate $\CIR$ without such slicing; however, the structure of $\vmu (y)$ makes this difficult. Recall that
\begin{equation} \label{eq:cond_exp2}
\vmu (y)
\;=\;
\int \vx \, d \sigma_{\vx|y} (\vx)
\end{equation}
is the conditional expectation of the inverse image of $f(\vx)$ for a fixed value of $y$. Approximating $\vmu (y)$ requires knowledge of the conditional measure $\sigma_{\vx|y}$, which may not be available if $f$ is complex. However, using the tools from Section \ref{sec:num_tools}, we can exploit composite structure in $\vmu (y)$.

The conditional expectation \eqref{eq:cond_exp2} is a function that maps values of $y$ to values in $\mathbb{R}^m$. Furthermore, $y$ is itself a function of $\vx$ (see \eqref{eq:y_eq_fx2}). Thus, the conditional expectation has composite structure. That is, we can define the function
\begin{equation}
\vmu_{\vx} (\vx)
\;=\;
\vmu (f(\vx)) ,
\end{equation}
where (using the notation from \eqref{eq:sX_sF_sG})
\begin{equation}
\begin{aligned}
f&: \left( \sX \subseteq \mathbb{R}^m \right) \rightarrow \left( \sF \subseteq \mathbb{R} \right) , \\
\vmu&: \left( \sF \subseteq \mathbb{R} \right) \rightarrow \left( \sG \subseteq \mathbb{R}^m \right) .
\end{aligned}
\end{equation}
Using the techniques from Section \ref{subsec:comp_func}, we look to construct a pseudospectral expansion of $\vmu$ using the orthogonal polynomials and Gauss-Christoffel quadrature with respect to $\gamma$.

Recall Assumption \ref{assume:stand_inputs} ensures the inputs have finite fourth moments and are standardized according to \eqref{eq:stand_inputs}. By Jensen's inequality,
\begin{equation}
\label{eq:Jensen_proof_LSIR1}
\vmu (y)^T \vmu (y)
\;\leq\;
\int \vx^T \vx \, d \sigma_{\vx|y} (\vx) , 
\end{equation}
Integrating both sides of \eqref{eq:Jensen_proof_LSIR1} with respect to $\gamma$,
\begin{equation}
\label{eq:Jensen_proof_LSIR2}
\int \vmu (y)^T \vmu (y) \, d \gamma(y)
\;\leq\;
\iint \vx^T \vx \, d \sigma_{\vx|y} (\vx)  \, d \gamma(y)
\;=\; \int \vx^T \vx \, d \rho (\vx) . 
\end{equation}
Under Assumption \ref{assume:stand_inputs}, the right-hand side of \eqref{eq:Jensen_proof_LSIR2} is finite. This guarantees that each component of $\vmu (y)$ is square integrable with respect to $\gamma$, ensuring that $\vmu (y)$ has a Fourier expansion in orthonormal polynomials with respect to $\gamma$,
\begin{equation} \label{eq:mu_expansion}
\vmu (y)
\;=\;
\sum_{i=0}^\infty \vmu_i \, \phi_i (y)
\end{equation}
with equality defined in the $L_2 (\gamma)$ sense. The Fourier coefficients in \eqref{eq:mu_expansion} are
\begin{equation} \label{eq:fourier_coeffs2}
\vmu_i
\;=\;
\int \vmu (y) \, \phi_i (y) \, d \gamma (y) .
\end{equation}
Plugging \eqref{eq:mu_expansion} into \eqref{eq:CIR2},
\begin{equation}
\begin{aligned}
\CIR \;&=\; \int \vmu (y) \, \vmu (y)^T \, d \gamma (y) \\
\;&=\; \int \left[ \sum_{i=0}^\infty \vmu_i \, \phi_i (y) \right] \left[ \sum_{j=0}^\infty \vmu_j \, \phi_j (y) \right]^T \, d \gamma (y) \\
\;&=\; \sum_{i=0}^\infty \sum_{j=0}^\infty \vmu_i \, \vmu_j^T \left[ \int \phi_i (y) \, \phi_j (y) \, d \gamma (y) \right] \\
\;&=\; \sum_{i=0}^\infty \vmu_i \, \vmu_i^T .
\end{aligned}
\end{equation}
Thus, the $\CIR$ matrix can be computed as the sum of the outer products of Fourier coefficients from \eqref{eq:fourier_coeffs2}.

We cannot compute the Fourier coefficients directly since they require knowledge of $\vmu (y)$. However, we can rewrite \eqref{eq:fourier_coeffs2} as
\begin{equation}
\begin{aligned}
\vmu_i
\;&=\;
\int \vmu (y) \, \phi_i (y) \, d \gamma (y) \\
\;&=\;
\int \left[ \int \vx \, d \sigma_{\vx|y} \right] \, \phi_i (y) \, d \gamma (y) \\
\;&=\;
\iint \vx \, \phi_i (y) \, d \sigma_{\vx|y} (\vx) \, d \gamma (y) \\
\;&=\;
\int \vx \, \phi_i (f(\vx)) \, d \rho (\vx) .
\end{aligned}
\end{equation}
This form of the Fourier coefficients is more amenable to numerical approximation. Since $\rho$ is known, we can obtain the $N$-point Gauss-Christoffel quadrature nodes $\vx_j$ and weights $\nu_j$ for $j = 0,\dots,N-1$, and define the pseudospectral coefficients
\begin{equation} \label{eq:quad_mu}
\hat{\vmu}_i
\;=\;
\sum_{j=0}^{N-1} \nu_j \, \vx_j \, \phi_i (f(\vx_j)) .
\end{equation}
To evaluate $\phi_i$ at $f (\vx_j)$, we use the Lanczos vectors in $\mV$. Recall from \eqref{eq:lanczos_vectors} that these vectors approximate the orthonormal polynomials from Stieltjes at $f$ evaluated at the quadrature nodes scaled by the square root of the associated quadrature weights. 

Algorithm \ref{alg:LSIR} provides an outline for Lanczos-Stieltjes inverse regression (LSIR). Note that this algorithm references the Lanczos algorithm for approximating composite functions (Algorithm \ref{alg:Lanc_comp_funcs}) which is defined for scalar-valued inputs. As discussed in Section \ref{subsec:comp_func}, we extend this to multivariate inputs by using tensor product constructions of the Gauss-Christoffel quadrature rule.

\begin{algorithm}[!ht]
\caption{Lanczos-Stieltjes inverse regression (LSIR)} \label{alg:LSIR}
\textbf{Given:} function $f: \mathbb{R}^m \rightarrow \mathbb{R}$ and input probability measure $\rho$ \\
\textbf{Assumptions:} Assumption \ref{assume:stand_inputs} holds \vspace{-0.8em}
\begin{enumerate}
\item Perform Algorithm \ref{alg:Lanc_comp_funcs} using $f$ and $\rho$ to obtain
\begin{equation}
\mA \mV \;=\; \mV \mT + \eta_k \vv_k \ve_k^T .
\end{equation} \vspace{-2.5em}
\item For $i = 0,\dots,m-1$, $\ell = 0,\dots,k-1$, \\
Compute the $i$th component of the $\ell$th pseudospectral coefficient
\begin{equation}
\left( \hat{\vmu}_\ell \right)_i
\;=\;
\sum_{p=0}^{N-1} \sqrt{\nu_p} \left( \vx_p \right)_i \left( \mV \right)_{p, \ell} .
\end{equation} \vspace{-2.5em}
\item For $i, j = 0,\dots,m-1$,\\
Compute the $i,j$th component of $\hCIR$
\begin{equation}
\label{eq:hCIR_inAlg}
\left( \hCIR \right)_{i,j}
\;=\;
\sum_{\ell=0}^{k-1} \left( \hat{\vmu}_\ell \right)_i \left( \hat{\vmu}_\ell \right)_j .
\end{equation} \vspace{-2em}
\end{enumerate}
\end{algorithm} 

Algorithm \ref{alg:LSIR} depends on two levels of approximation: (i) approximation due to the quadrature rule over $\sX$ and (ii) approximation due to truncating the polynomial expansion of $\vmu (y)$. The former depends on the number $N$ of nodes used while the latter depends on the number $k$ of Lanczos iterations performed. Performing more Lanczos iterations includes more terms in the approximation of $\CIR$; however, the additional terms correspond to integrals against higher degree polynomials in \eqref{eq:fourier_coeffs2}. For a fixed number of quadrature nodes, the quality of the pseudospectral approximation deteriorates as the degree of polynomial increases. Therefore, sufficiently many quadrature nodes are needed to ensure quality estimates of the orthonormal polynomials of high degree. 


\subsection{Lanczos-Stieltjes average variance estimation (LSAVE)}
\label{subsec:LSAVE}

In this section, we apply the Lanczos-Stieltjes approach from Section \ref{subsec:LSIR} to the $\CAVE$ matrix to construct an alternative algorithm to the slice-based SAVE. Recall from \eqref{eq:CIR_CAVE} that
\begin{equation} \label{eq:CAVE2}
\CAVE
\;=\;
\int \left( \mI - \mSigma (y) \right)^2 \, d \gamma (y) ,
\end{equation}
where
\begin{equation}
\label{eq:cond_cov2}
\mSigma (y)
\;=\;
\int \left( \vx - \vmu (y) \right) \left( \vx - \vmu (y) \right)^T \, d \sigma_{\vx|y} (\vx)
\end{equation}
is the conditional covariance of the inverse image of $f(\vx)$ for a fixed value of $y$. Similar to the conditional expectation, $\mSigma (y)$ has composite structure due to the relationship $y = f(\vx)$ such that we can define
\begin{equation}
\mSigma_{\vx} (\vx)
\;=\;
\mSigma (f(\vx))  ,
\end{equation}
where
\begin{equation}
\begin{aligned}
f&: \left( \sX \subseteq \mathbb{R}^m \right) \rightarrow \left( \sF \subseteq \mathbb{R} \right) , \\
\mSigma&: \left( \sF \subseteq \mathbb{R} \right) \rightarrow \left( \sG \subseteq \mathbb{R}^{m \times m} \right) .
\end{aligned}
\end{equation}
We want to build a pseudospectral expansion of $\mSigma$ with respect to $\gamma$ similar to $\vmu$ in Section \ref{subsec:LSIR}.

Recall Assumption \ref{assume:stand_inputs} as we consider the Frobenius norm of the conditional covariance. By Jensen's inequality,
\begin{equation}
\left| \left| \mSigma (y) \right| \right|_F^2
\;\leq\;
\int \left| \left| \left( \vx - \vmu(y) \right) \left( \vx - \vmu(y) \right)^T \right| \right|_F^2 \, d \sigma_{\vx|y} (\vx) .
\end{equation}
Integrating each side with respect to $\gamma$,
\begin{equation}
\begin{aligned}
\int \left| \left| \mSigma (y) \right| \right|_F^2 \, d \gamma (y)
\;&\leq\;
\iint \left| \left| \left( \vx - \vmu(y) \right) \left( \vx - \vmu(y) \right)^T \right| \right|_F^2 \, d \sigma_{\vx|y} (\vx) \, d \gamma (y) \\
\;&\leq\;
\int \left| \left| \left( \vx - \vmu(f(\vx)) \right) \left( \vx - \vmu(f(\vx)) \right)^T \right| \right|_F^2 \, d \rho (\vx) .
\end{aligned}
\end{equation}
Expanding the integrand on the right-hand side produces sums and products of fourth and lower conditional moments of the inputs. Assumption \ref{assume:stand_inputs} guarantees that all of these conditional moments are finite. Thus,
\begin{equation}
\int \left| \left| \mSigma (y) \right| \right|_F^2 \, d \gamma (y)
\;<\;
\infty ,
\end{equation}
which implies that each component of $\mSigma (y)$ is square integrable with respect to $\gamma$; therefore it has a convergent Fourier expansion in terms of orthonormal polynomials with respect to $\gamma$, 
\begin{equation} \label{eq:sigma_expansion}
\mSigma (y)
\;=\;
\sum_{i=0}^\infty \mSigma_i \, \phi_i (y) ,
\end{equation}
where equality is in the $L_2 (\gamma)$ sense. The coefficients in \eqref{eq:sigma_expansion} are
\begin{equation}
\label{eq:Sigma_fourier_coeff}
\mSigma_i
\;=\;
\int \mSigma (y) \, \phi_i (y) \, d \gamma (y) .
\end{equation}
Plugging \eqref{eq:sigma_expansion} into \eqref{eq:CAVE2}
\begin{equation}
\begin{aligned}
\CAVE
\;&=\;
\int \left( \mI - \mSigma (y) \right)^2 \, d \gamma (y) \\
\;&=\;
\int \left( \mI - \sum_{i=0}^\infty \mSigma_i \, \phi_i (y) \right)^2 \, d \gamma (y) \\
\;&=\;
\int \mI \, d \gamma (y) - 2 \sum_{i=0}^{\infty} \mSigma_i \int \phi_i (y) \, d \gamma (y) + \sum_{i=0}^\infty \sum_{j=0}^\infty \mSigma_i \, \mSigma_j \int \phi_i (y) \, \phi_j (y) \, d \gamma (y) \\
\;&=\;
\mI - 2 \, \mSigma_0 + \sum_{i=0}^\infty \mSigma_i^2 .
\end{aligned}
\end{equation}
Therefore, we can compute $\CAVE$ using the Fourier coefficients of $\mSigma (y)$. To simplify the computation of $\mSigma_i$, we rewrite \ref{eq:Sigma_fourier_coeff} as
\begin{equation}
\begin{aligned}
\mSigma_i
\;&=\;
\int \mSigma (y) \, \phi_i (y) \, d \gamma (y) \\
\;&=\;
\int \left( \int \left( \vx - \vmu (y) \right) \left( \vx - \vmu (y) \right)^T \, d \sigma_{\vx|y} (\vx) \right) \phi_i (y) \, d \gamma (y) \\
\;&=\;
\int \int \left( \vx - \vmu (y) \right) \left( \vx - \vmu (y) \right)^T \phi_i (y) \, d \sigma_{\vx|y} (\vx) \, d \gamma (y) \\
\;&=\;
\int \left( \vx - \vmu (f(\vx)) \right) \left( \vx - \vmu (f(\vx)) \right)^T \phi_i (f(\vx)) \, d \rho (\vx) .
\end{aligned}
\end{equation}
We approximate this integral using the $N$-point Gauss-Christoffel quadrature rule with respect to $\rho$ to obtain
\begin{equation} \label{eq:quad_Sigma}
\hat{\mSigma}_i
\;=\;
\sum_{j=0}^{N-1} \nu_j \left( \vx_j - \vmu (f(\vx_j)) \right) \left( \vx_j - \vmu (f(\vx_j)) \right)^T \phi_i (f(\vx_j)) .
\end{equation}
We again approximate $\phi_i (f(\vx_j))$ using the Lanczos vectors similar to LSIR. Notice that \eqref{eq:quad_Sigma} also depends on $\vmu (f(\vx_j))$ for $j = 0, \dots ,N-1$. To obtain these values, we compute the pseudospectral coefficients of $\vmu (f(\vx))$ from \eqref{eq:quad_mu} and construct its pseudospectral expansion at each $\vx_j$. Algorithm \ref{alg:LSAVE} provides an outline for Lanczos-Stieltjes average variance estimation (LSAVE).

\begin{algorithm}[!ht]
\caption{Lanczos-Stieltjes average variance estimation (LSAVE)} \label{alg:LSAVE}
\textbf{Given:} function $f: \mathbb{R}^m \rightarrow \mathbb{R}$ and input probability measure $\rho$ \\
\textbf{Assumptions:} Assumption \ref{assume:stand_inputs} holds \vspace{-0.8em}
\begin{enumerate}
\item Perform Algorithm \ref{alg:Lanc_comp_funcs} using $f$ and $\rho$ to obtain
\begin{equation}
\mA \mV \;=\; \mV \mT + \eta_k \vv_k \ve_k^T .
\end{equation} \vspace{-2.5em}
\item For $i = 0,\dots,m-1$, $\ell = 0,\dots,k-1$, \\
Compute the $i$th component of the $\ell$th pseudospectral coefficient of $\vmu (y)$
\begin{equation}
\left( \hat{\vmu}_\ell \right)_i
\;=\;
\sum_{p=0}^{N-1} \sqrt{\nu_p} \left( \vx_p \right)_i \left( \mV \right)_{p, \ell}.
\end{equation} \vspace*{-2.5em}
\item For $i = 0,\dots,m-1$, \\
Compute the $i$th component of the pseudospectral expansion of $\vmu (f(\vx_p))$
\begin{equation}
\left( \hat{\vmu} (f (\vx_p)) \right)_i
\;=\;
\sum_{\ell = 0}^{k-1} \frac{1}{\sqrt{\nu_p}} \, \left( \hat{\vmu}_\ell \right)_i \left( \mV \right)_{p, \ell} .
\end{equation} \vspace*{-2.5em}
\item For $i,j = 0,\dots,m-1$, $\ell = 0,\dots,k-1$, \\
Compute the $i,j$th component of the $\ell$th pseudospectral coefficient of $\mSigma (y)$
\begin{equation}
\left( \hat{\mSigma}_\ell \right)_{i,j}
\;=\;
\sum_{p=0}^{N-1} \sqrt{\nu_p} \left( (\vx_p)_i - \left( \hat{\vmu} (f (\vx_p)) \right)_i \right) \left( (\vx_p)_j - \left( \hat{\vmu} (f (\vx_p)) \right)_j \right) \left( \mV \right)_{p, \ell} .
\end{equation} \vspace*{-2.5em}
\item For $i,j = 0,\dots,m-1$, \\
Compute the $i,j$th component of $\hCAVE$
\begin{equation}
\left( \hCAVE \right)_{i,j}
\;=\;
\delta_{i,j} - 2 \, \left( \hat{\mSigma}_0 \right)_{i,j} + \sum_{\ell=0}^{k-1} \sum_{p=0}^{m-1} \left( \hat{\mSigma}_\ell \right)_{i,p} \left( \hat{\mSigma}_\ell \right)_{p,j} .
\end{equation} \vspace*{-2em}
\end{enumerate}
\end{algorithm}

Note that Algorithm \ref{alg:LSAVE} contains the same two-level approximation as in the LSIR algorithm. As such, it also requires a sufficiently high-order quadrature rule to accurately approximate the high degree polynomials resulting from $k$ Lanczos iterations. In the next section, we provide numerical studies of the LSIR and LSAVE algorithms on several test problems as well as comparisons to the traditional SIR and SAVE algorithms.

%% file: sec6-num_results.tex
\section{Numerical results}
\label{sec:numericalresults}


In this section, we numerically study the LSIR and LSAVE algorithms. In Section \ref{subsec:LS_converge}, we study the approximation of the Gauss-Christoffel quadrature and orthonormal polynomials by the Lanczos algorithm as described in Section \ref{subsec:comp_func}. This study is performed on a multivariate quadratic function of dimension $m = 3$. In Sections \ref{subsec:1d_ridge} and \ref{subsec:OTL_circuit}, we examine the approximation of $\CIR$ and $\CAVE$ by the Lanczos-Stieltjes estimates $\hCIR$ and $\hCAVE$, and we compare these results to the slice-based estimates $\hCSIR$ and $\hCSAVE$. This study is performed on two problems: (i) an exact one-dimensional ridge function with $m = 5$ (see Section \ref{subsec:1d_ridge}) and (ii) a model of the output from a transformerless push-pull circuit with $m = 6$ inputs (see Section \ref{subsec:OTL_circuit}). Matlab code for performing these studies is available at \url{https://bitbucket.org/aglaws/gauss-christoffel-quadrature-for-inverse-regression}

Recall that the goal of the proposed algorithms is to approximate the $\CIR$ and $\CAVE$ matrices better than the slice-based (Riemann sum) approach in SIR and SAVE. Therefore, most of the presented results are in terms of the relative matrix errors and differences in the Frobenius norm. Recall that the Frobenius norm is given by
\begin{equation}
\left| \left| \mE \right| \right|_F
\;=\;
\left( \sum_{i=0}^{m-1} \sum_{j=0}^{m-1} \left( \mE \right)_{i,j}^2 \right)^{1/2} .
\end{equation}
References to matrix differences or matrix errors in the following examples should be interpreted to mean the Frobenius norm of these differences or errors.


\subsection{Example 1: Lanczos-Stieltjes convergence}
\label{subsec:LS_converge}

In this section, we consider
\begin{equation}
\label{eq:ex1}
y
\;=\;
f(\vx)
\;=\;
\vg^T \vx + \vx^T \mH \vx ,
\qquad
\vx \in [-1, 1]^3 \subset \mathbb{R}^3 ,
\end{equation}
where $\vg \in \mathbb{R}^3$ is a constant vector and $\mH \in \mathbb{R}^{3 \times 3}$ is a constant matrix. We assume the inputs are weighted by the uniform density over the input space $\sX = [-1, 1]^3$,
\begin{equation}
d \rho (\vx)
\;=\;
\begin{cases} 
\frac{1}{2^3} \, d \vx & \text{if } ||\vx||_\infty \leq 1 , \\
0 \, d \vx & \text{otherwise.}
\end{cases}
\end{equation} 

Recall from Section \ref{subsec:comp_func} that Algorithm \ref{alg:Lanc_comp_funcs} produces a $k$-point Gauss-Christoffel quadrature rule and the first $k$ orthonormal polynomials relative to the discrete measure $\gamma_N$. We treat these as approximations to the quadrature rule and orthonormal polynomials relative to the continuous measure $\gamma$ (see \eqref{eq:approx_quad} and \eqref{eq:lanczos_vectors}). In this section, we study the behavior of these approximations for \eqref{eq:ex1}.

\begin{figure}[!ht]
\centering
\includegraphics[width=1\textwidth]{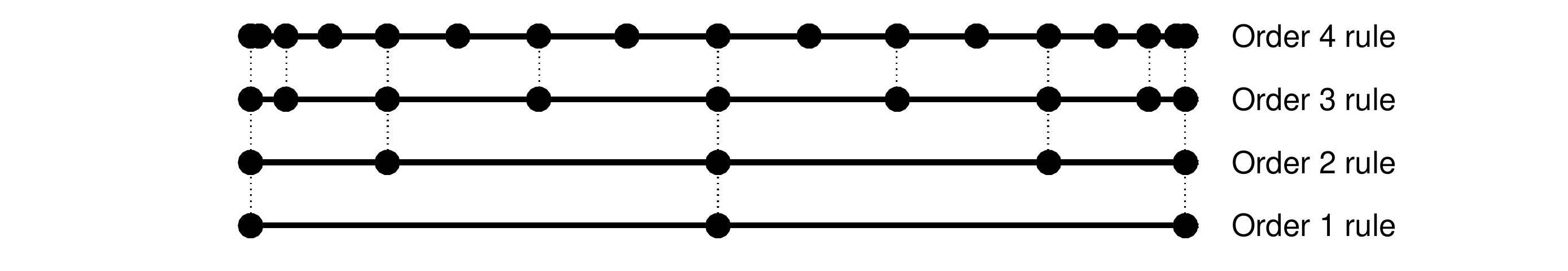}
\caption{A visual representation of the (one-dimensional) nested Clenshaw-Curtis quadrature rule. The order $M$ rule contains all the nodes from the order $M-1$ rule. The studies in this section use a tensor product Clenshaw-Curtis quadrature rule.}
\label{fig:nested_quad}
\end{figure}

We use a nested quadrature rule for this study. An order $M$ nested quadrature rule contains all the nodes used in the order $M-1$ rule (see Figure \ref{fig:nested_quad}). This enables direct comparison in studying the Lanczos vectors, which approximate evaluations of the orthonormal polynomials at the $f_i$'s associated with the input quadrature rule. We use a tensor product Clenshaw-Curtis quadrature rule on $\sX$ for this study~\citep{Clenshaw60}. Note that the Gauss-Christoffel quadrature rule used elsewhere in this paper is not in general nested. 

\begin{figure}[!ht]
\centering
\subfloat[]{
\label{fig:lambda_convergence}
\includegraphics[width=0.45\textwidth]{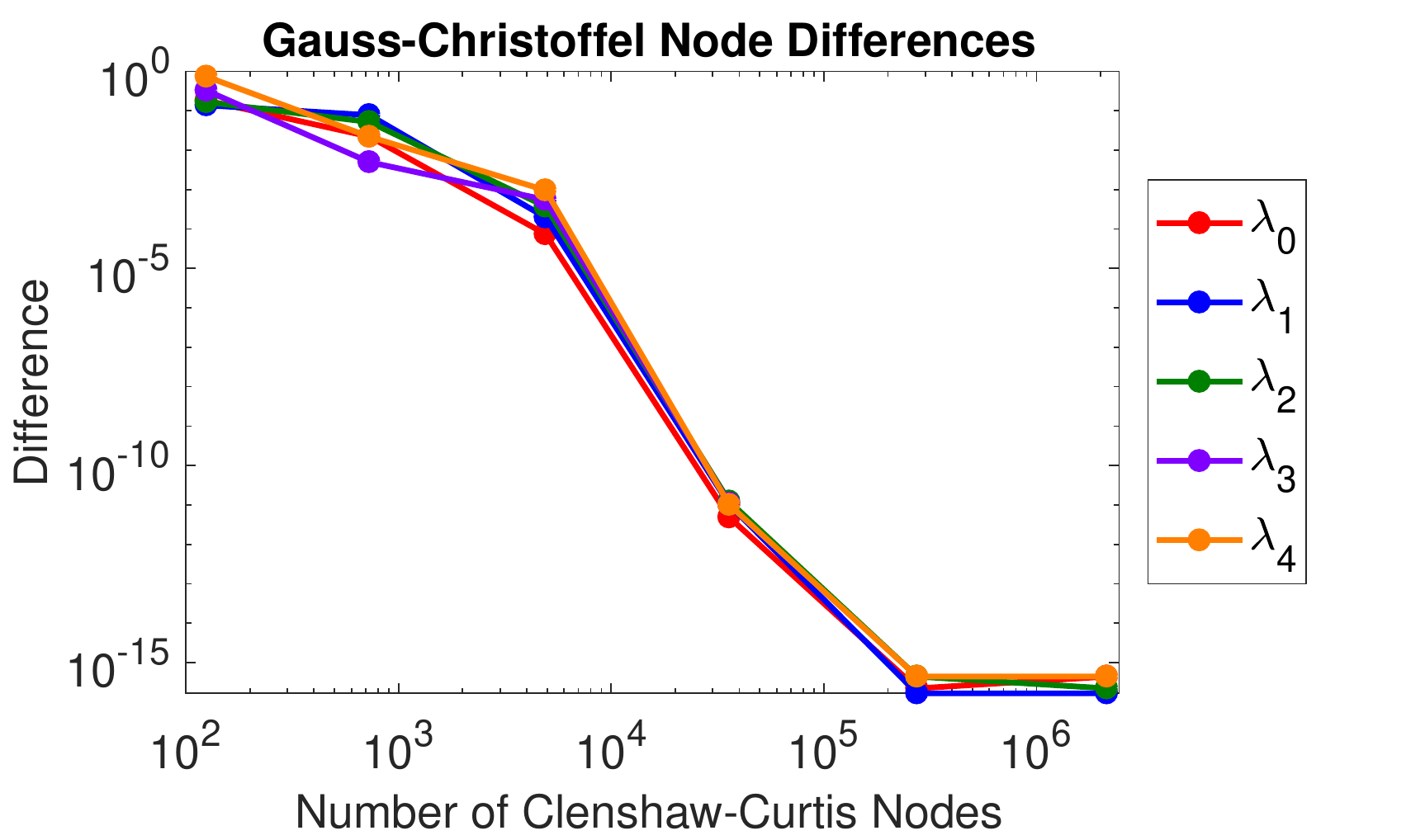}
}
\subfloat[]{
\label{fig:omega_convergence}
\includegraphics[width=0.45\textwidth]{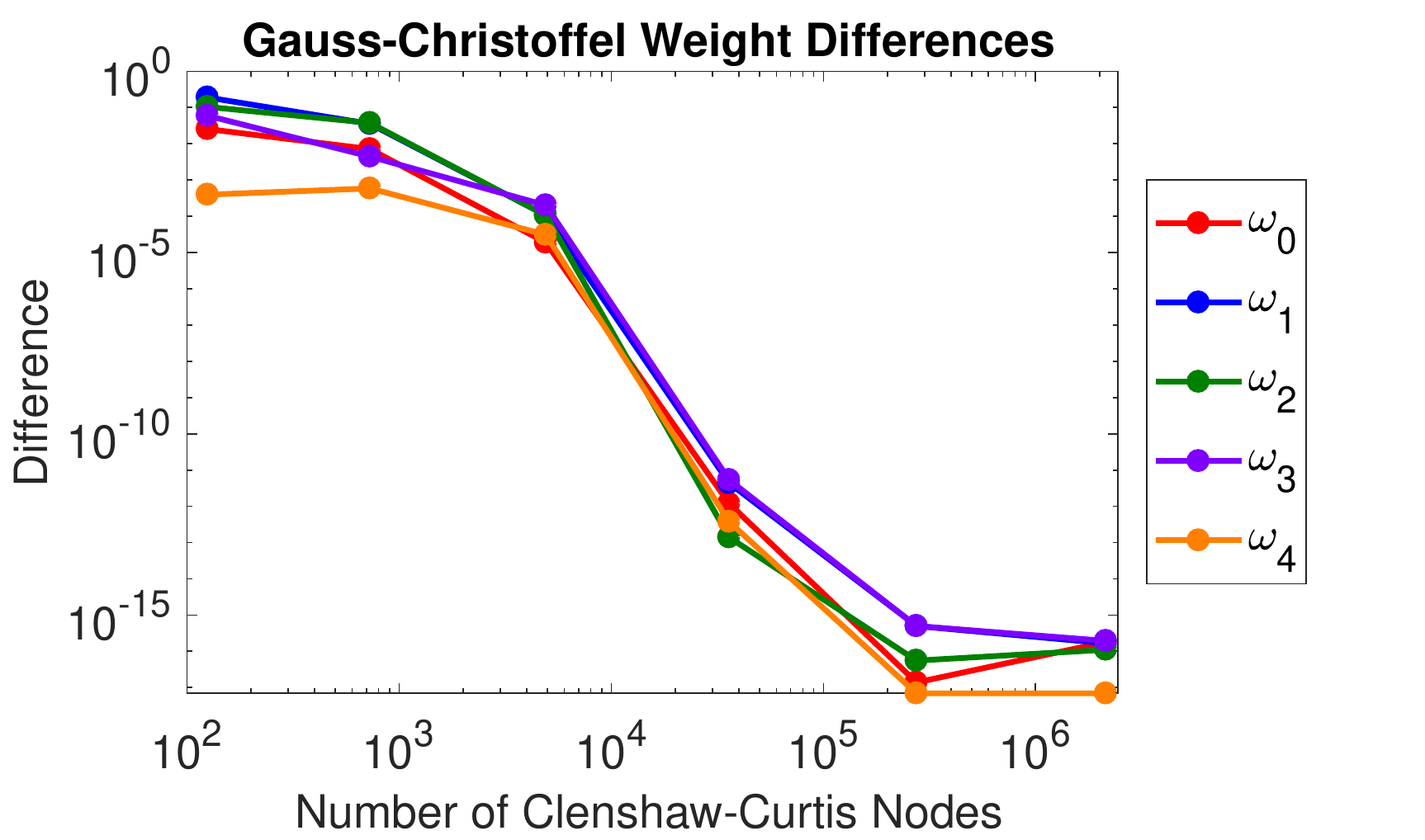}
}
\caption{Subsequent differences in the 5-point Gauss-Christoffel quadrature nodes (Figure \ref{fig:lambda_convergence}) and weights (Figure \ref{fig:omega_convergence}) resulting from Algorithm \ref{alg:Lanc_comp_funcs} applied to \eqref{eq:ex1} with increasing orders of tensor product Clenshaw-Curtis quadrature rules used over $\sX$. As the orders increase, the approximate quadrature rule over $\sF$ converges.}
\label{fig:quad_convergence}
\end{figure}

First, we study convergence of the Gauss-Christoffel quadrature rule produced on the output space $\sF$ by the Lanczos algorithm. \citet{Gautschi04} shows convergence of the Jacobi matrix \eqref{eq:T_Jacobi} to \eqref{eq:J_Jacobi} as the discrete approximation $\gamma_N$ approaches $\gamma$; however, we are specifically interested in the quadrature rule resulting from the eigendecomposition of the Jacobi matrix. Let $\hat{\lambda}_i^{(M)}$, $\hat{\omega}_i^{(M)}$ denote the $i$th estimated Gauss-Christoffel quadrature node and weight resulting from Algorithm \ref{alg:Lanc_comp_funcs} with an order $M$ Clenshaw-Curtis quadrature on the input space. We are interested in the differences in the $k$-point Gauss-Christoffel quadrature nodes and weights resulting from subsequent orders of Clenshaw-Curtis quadrature on the input space $\sX$,
\begin{equation}
\left| \hat{\lambda}_i^{(M)} - \hat{\lambda}_i^{(M-1)} \right| \quad \text{and} \quad \left| \hat{\omega}_i^{(M)} - \hat{\omega}_i^{(M-1)} \right|,
\end{equation}
for $i = 0,\dots,k-1$. Figure \ref{fig:quad_convergence} contains plots of these differences for $k = 5$. The Gauss-Christoffel quadrature rule with respect to $\gamma$ converges as higher-order Clenshaw-Curtis quadrature rules are used on the input space. Quality estimates of the quadrature rule over the output space $\sF$ are required to produce good approximations of the $\CIR$ and $\CAVE$ using the Lanczos-Stieltjes approach.

\begin{figure}[!ht]
\centering
\includegraphics[width=0.45\textwidth]{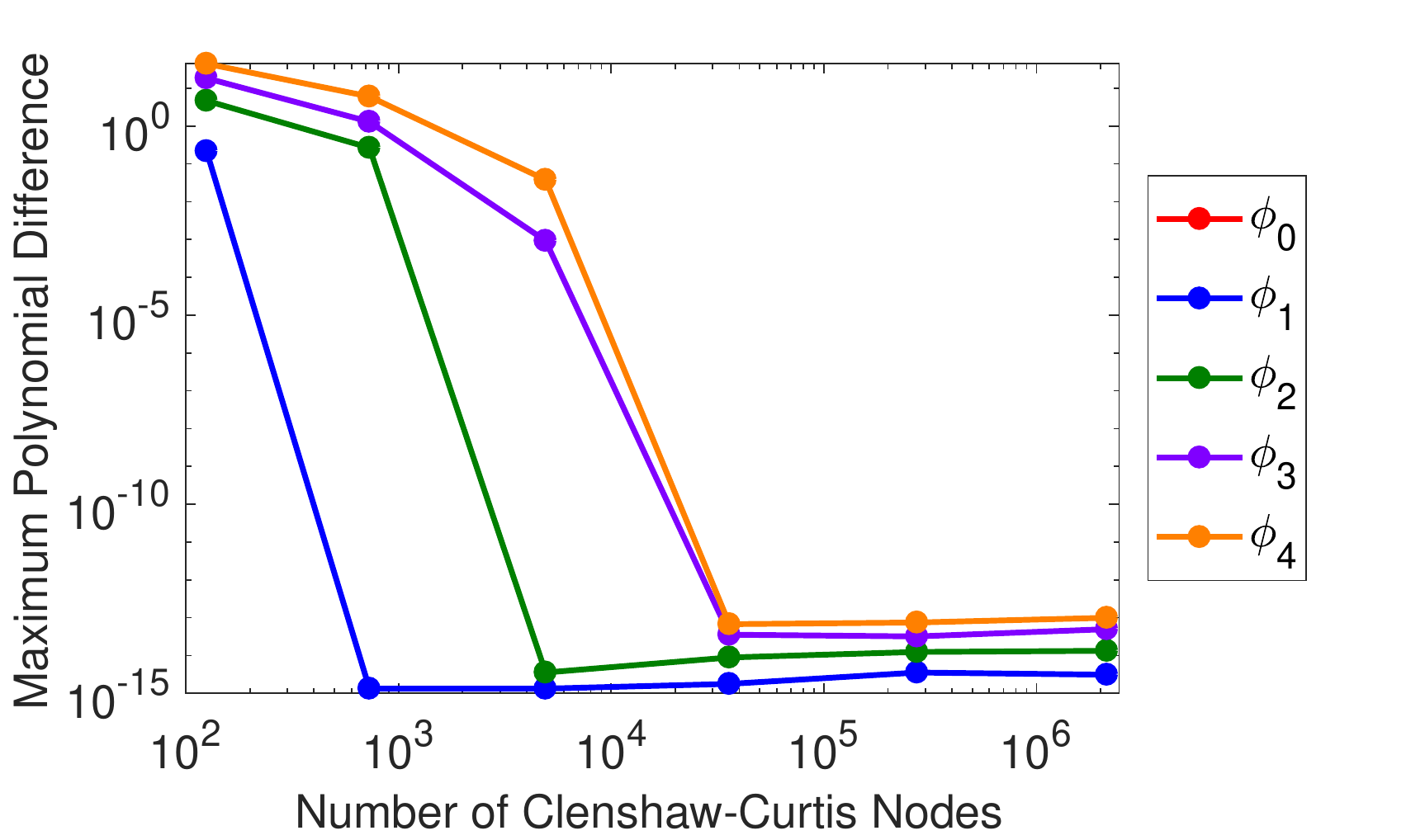}
\caption{Maximum differences in the first 5 approximated orthonormal polynomials as we increase the order of the quadrature rule used to approximate them. Higher degree polynomials require more quadrature nodes to be accurately estimated than lower degree polynomials.}
\label{fig:phi_converge}
\end{figure}

Next, we examine convergence of the Lanczos vectors to the orthonormal polynomials with respect to $\gamma$. Recall from \eqref{eq:lanczos_vectors} that the Lanczos vectors contain evaluations of the first $k$ (corresponding to the number of Lanczos iterations performed) orthonormal polynomials with respect to $\gamma_N$ at the $f_i$'s from Algorithm \ref{alg:Lanc_comp_funcs} weighted by the square root of the associated quadrature weights, $\sqrt{\nu_i}$. Let $\mV^{(M)}$ be the matrix of $k$ Lanczos vectors resulting from the order $M$ Clenshaw-Curtis quadrature rule, and define $\mW_{\nu} = \diag{[ \, \sqrt{\nu_0} \, \dots \, \sqrt{\nu_{N-1}} \, ]}$. We define $\tilde{\mV}^{(M)} = \mW_{\nu}^{-1} \mV^{(M)}$ to be the matrix of orthonormal polynomials evaluated at the $f_i$'s (no longer scaled by the $\sqrt{\nu_i}$'s). This matrix is $N_M \times k$, where $N_M$ is the number of nodes resulting from the order $M$ tensor product Clenshaw-Curtis quadrature rule. The matrix $\tilde{\mV}^{(M-1)}$ contains evaluations of the same $k$ orthonormal polynomials at a subset of the quadrature nodes. Let $\mP \in \mathbb{R}^{N_{M-1} \times N_M}$ be the matrix which removes the rows of $\tilde{\mV}^{(M)}$ that do not correspond to rows in $\tilde{\mV}^{(M-1)}$. We write the maximum difference in the $i$th polynomial for subsequent orders of the quadrature rule as
\begin{equation}
\left| \left| \mP \left(\tilde{\mV}^{(M)}\right)_i - \left(\tilde{\mV}^{(M-1)}\right)_i \right| \right|_{\infty} , \quad i = 0, \dots, k-1 ,
\end{equation}
where $(\cdot)_i$ denotes the $i$th column of the given matrix. Figure \ref{fig:phi_converge} contains a plot of these maximum polynomial differences for increasing nested quadrature orders. We note that higher degree polynomials require a higher order rule to produce accurate approximations. This is not surprising, but it does highlight an important relationship between the number of quadrature nodes taken with respect to $\rho$ and the number of Lanczos iterations performed. Namely, as the number of Lanczos iterations increases, more quadrature nodes are required to ensure accurate approximation of the orthonormal polynomials (and, in turn, $\CIR$ and $\CAVE$).


\subsection{Example 2: One-dimension ridge function}
\label{subsec:1d_ridge}

In this section, we study the convergence of the LSIR and LSAVE methods and compare these methods to their slice-based counterparts. Consider inputs $\vx \in \mathbb{R}^5$ weighted by a standard Gaussian,
\begin{equation}
d \rho (\vx)
\;=\;
\frac{1}{(2 \pi)^{5/2}} e^{-(\vx^T \vx) / 2} d \vx ,
\end{equation}
and let
\begin{equation} \label{eq:ex2}
y
\;=\;
f(\vx)
\;=\;
\va^T \vx \, \cos \left( \frac{\va^T \vx}{2 \pi} \right) ,
\end{equation}
where $\va \in \mathbb{R}^5$ is a constant vector. It is clear from inspection that $f(\vx)$ is a one-dimensional ridge function in terms of $\va^T \vx \in \mathbb{R}$. 

\begin{figure}[!ht]
\centering
\subfloat[]{
\label{fig:LSIR_quadpts_Ex2}
\includegraphics[width=0.45\textwidth]{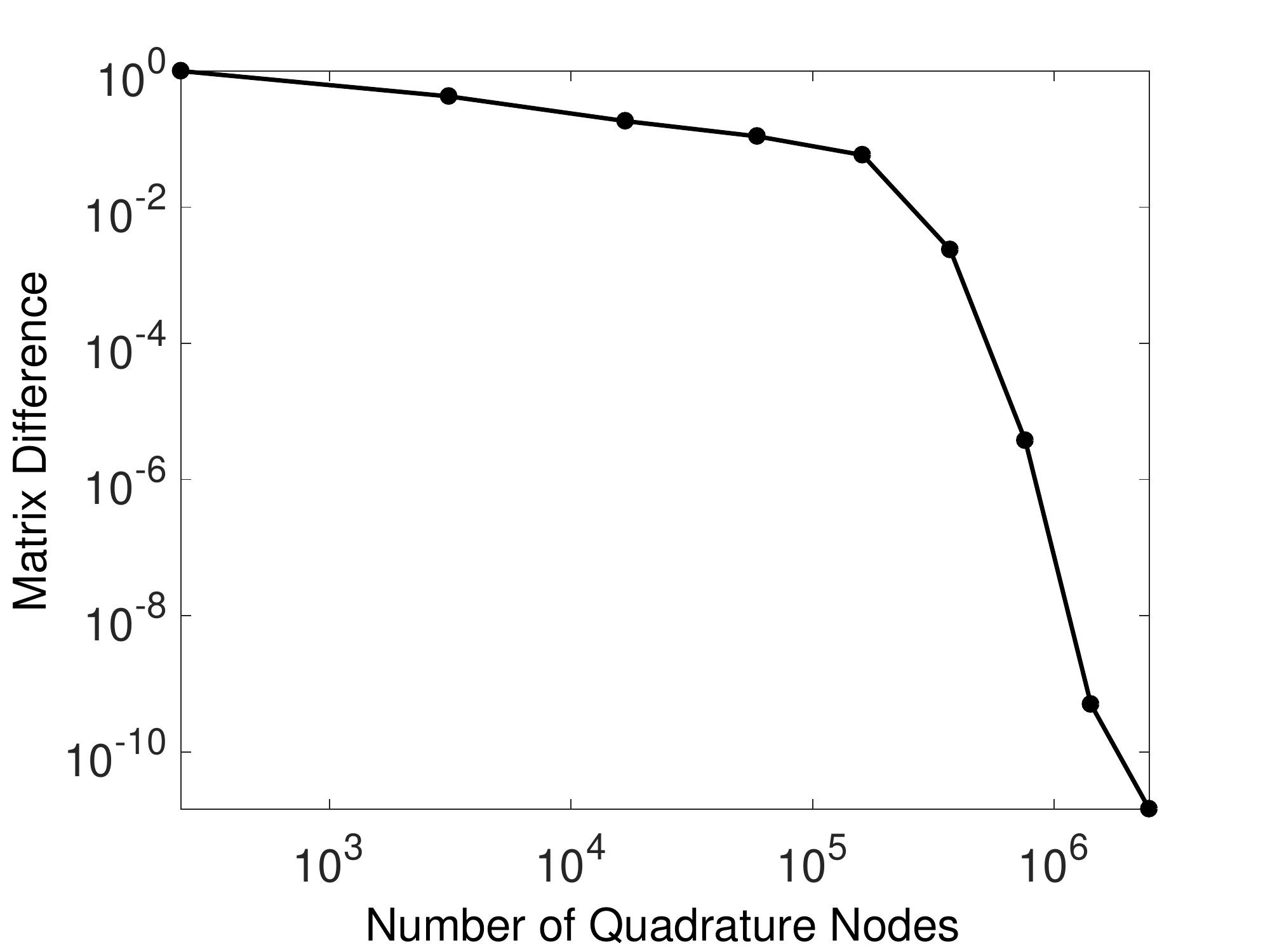}
}
\subfloat[]{
\label{fig:LSIR_lanciters_Ex2}
\includegraphics[width=0.45\textwidth]{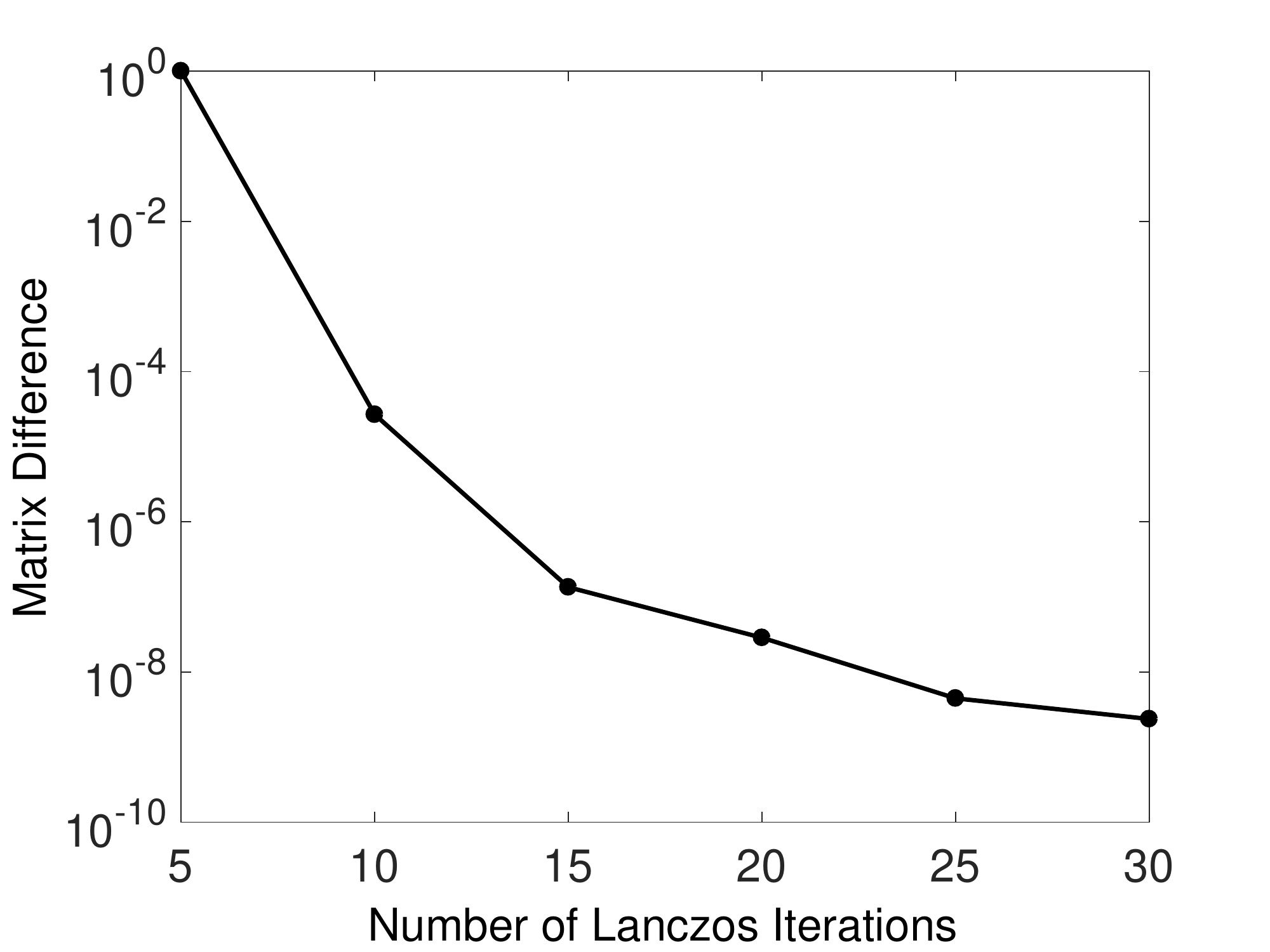}
}
\caption{Subsequent matrix differences for increasing numbers of quadrature nodes with $k = 35$ (Figure \ref{fig:LSIR_quadpts_Ex2}) and for increasing numbers of Lanczos iterations with $N = 21^5 =$ 4,084,101 quadrature nodes (Figure \ref{fig:LSIR_lanciters_Ex2}) on \eqref{eq:ex2}.}
\label{fig:LSIR_relconv_Ex2}
\end{figure}

We first study the convergence of the LSIR algorithm in terms of both the number of Gauss-Christoffel quadrature nodes over the input space ($N$ in Algorithm \ref{alg:Lanc_comp_funcs}) and the number of Lanczos iterations performed ($k$ in Algorithm \ref{alg:Lanc_comp_funcs}). Figure \ref{fig:LSIR_relconv_Ex2} contains the results of these studies. The first plot shows differences between subsequent Lanczos-Stieltjes approximations of $\CIR$ computed using increasing numbers of quadrature nodes with the number of Lanczos iterations fixed at $k = 35$. We see a decay in these differences as we increase the number of quadrature nodes. The second plot in Figure \ref{fig:LSIR_relconv_Ex2} contains subsequent matrix differences for increasing numbers of Lanczos iterations with the number of quadrature nodes fixed at 21 per dimension. Recall that we use a tensor product extension of the one-dimensional Gauss-Christoffel quadrature rule to $m = 5$ dimensions. This results in $N = 21^5 =$ 4,084,101 total quadrature nodes. We again see convergence in the differences as we increase $k$. Recall that the number of Lanczos iterations corresponds to the degree of the polynomial approximation of $\vmu (y)$ (see \eqref{eq:hCIR_inAlg}).

\begin{figure}[!ht]
\centering
\includegraphics[width=0.45\textwidth]{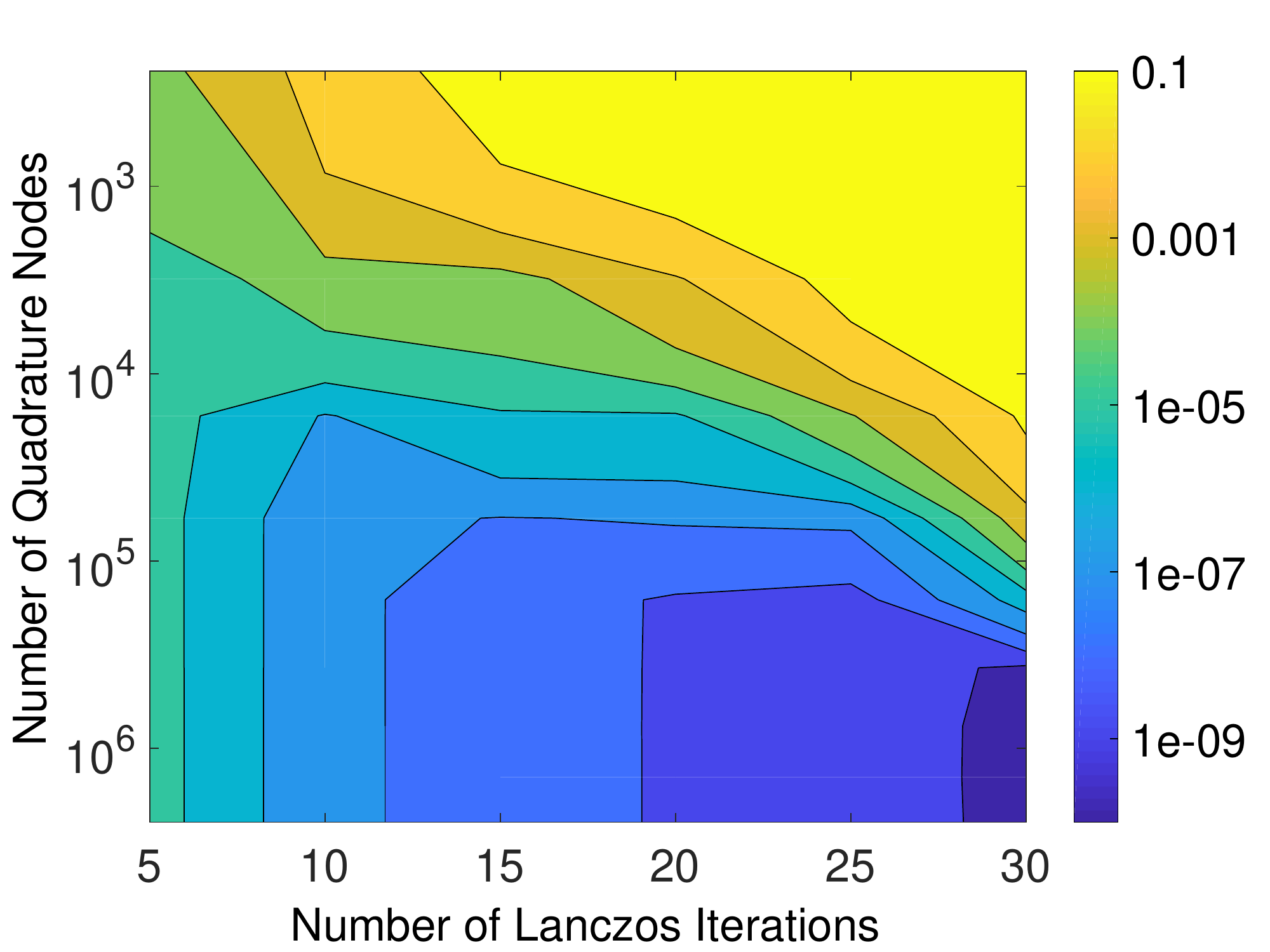}
\caption{Study of the matrix errors for various values of $N$ and $k$ relative to the ``true'' $\CIR$ matrix (i.e., $N = 21^5 =$ 4,084,101, $k = 35$) for \eqref{eq:ex2}.}
\label{fig:LSIR_2d_converge_Ex2}
\end{figure}

The results from Figure \ref{fig:LSIR_relconv_Ex2} suggest that the LSIR algorithm is converging in terms of $N$ and $k$. We use this as justification for the next study in which we perform Algorithm \ref{alg:LSIR} for $N = 21^5 =$ 4,084,101 and $k = 35$ and treat the resulting matrix as the ``true'' value of $\CIR$. We then compute errors relative to this matrix for various values of $N$ and $k$. The results of this study are contained in Figure \ref{fig:LSIR_2d_converge_Ex2}. 

Figure \ref{fig:LSIR_2d_converge_Ex2} shows the relative error decaying as we increase both the number of quadrature nodes and Lanczos iterations (i.e., as we move down and to the right). Consider this plot for increasing $N$ with $k$ fixed (i.e., moving downward at a fixed point along the horizontal axis). The error decays up to a point at which it remains constant. This decay corresponds to more accurate computation of the pseudospectral coefficients $\hat{\vmu}_i$ by taking more quadrature nodes. The leveling off corresponds to the point at which errors in the coefficients are smaller than errors due to truncating the pseudospectral expansion at $k$ (the number of Lanczos iterations). Conversely, if we fix $N$ and study the error as we increase $k$ (i.e., fix a point along the vertical axis and move right), we see the error decay until a point at which it begins to grow again. This behavior agrees with the results from Section \ref{subsec:LS_converge} which suggest that sufficiently many quadrature nodes are needed to accurately estimate the high-degree polynomials associated with large values of $k$. As we move right in Figure \ref{fig:LSIR_2d_converge_Ex2}, we are approximating higher-degree polynomials using the Lanczos algorithm. Poor approximations of these polynomials results in an inaccurate estimates of $\CIR$.

\begin{figure}[!ht]
\centering
\subfloat[]{
\label{fig:LSAVE_quadpts_Ex2}
\includegraphics[width=0.4\textwidth]{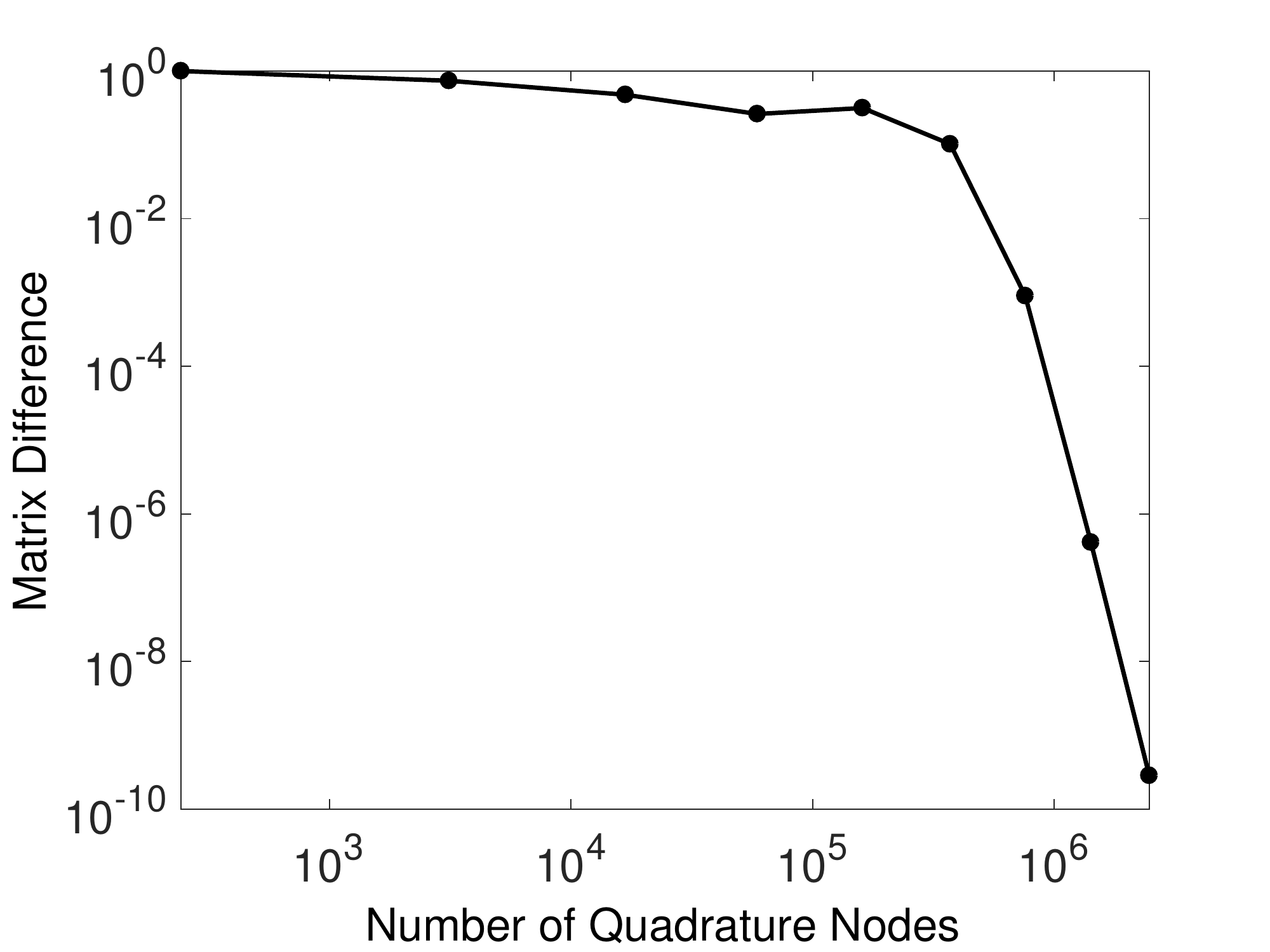}
}
\subfloat[]{
\label{fig:LSAVE_lanciters_Ex2}
\includegraphics[width=0.4\textwidth]{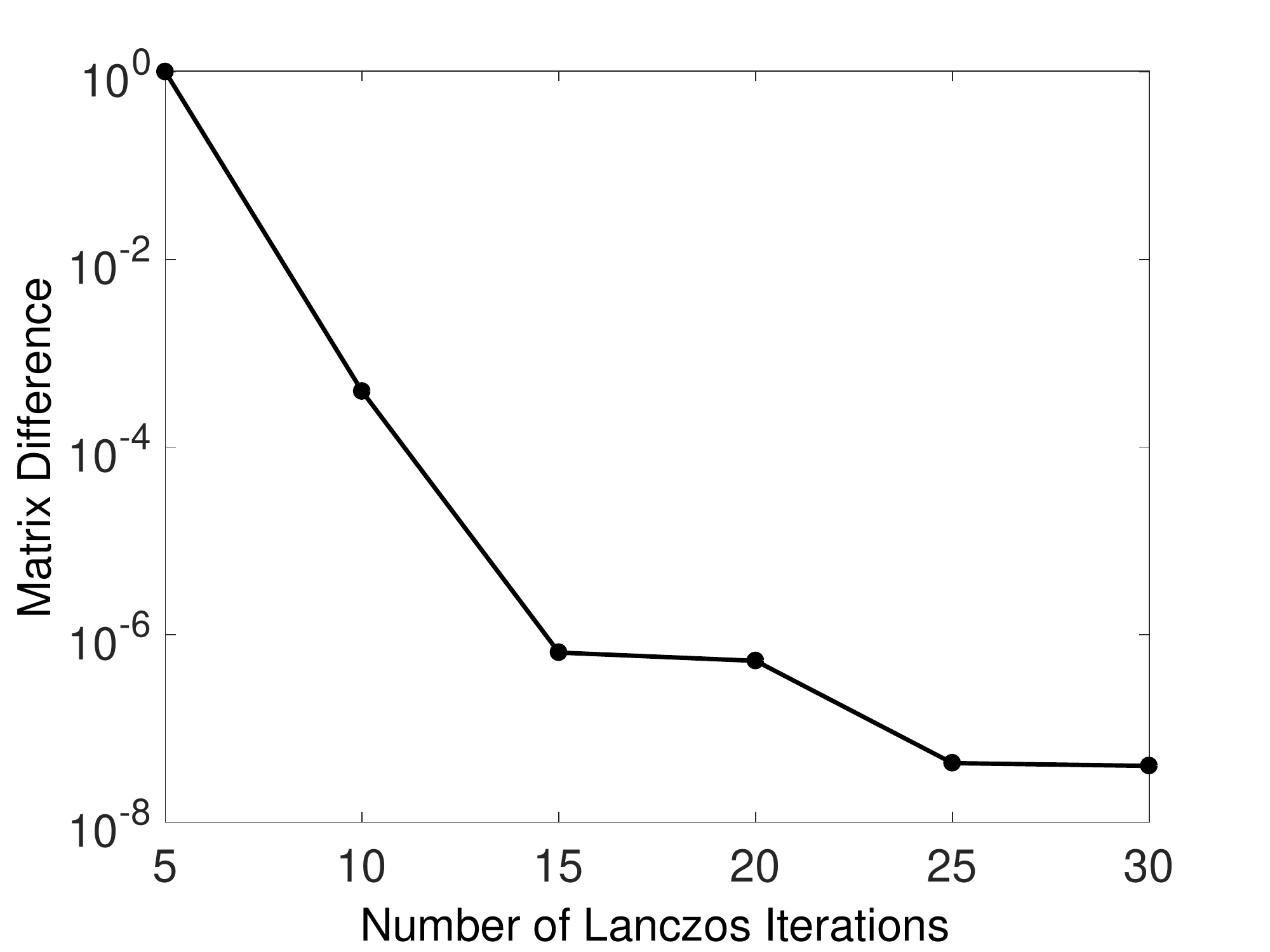}
} \\
\subfloat[]{
\label{fig:LSAVE_2d_converge_Ex2}
\includegraphics[width=0.4\textwidth]{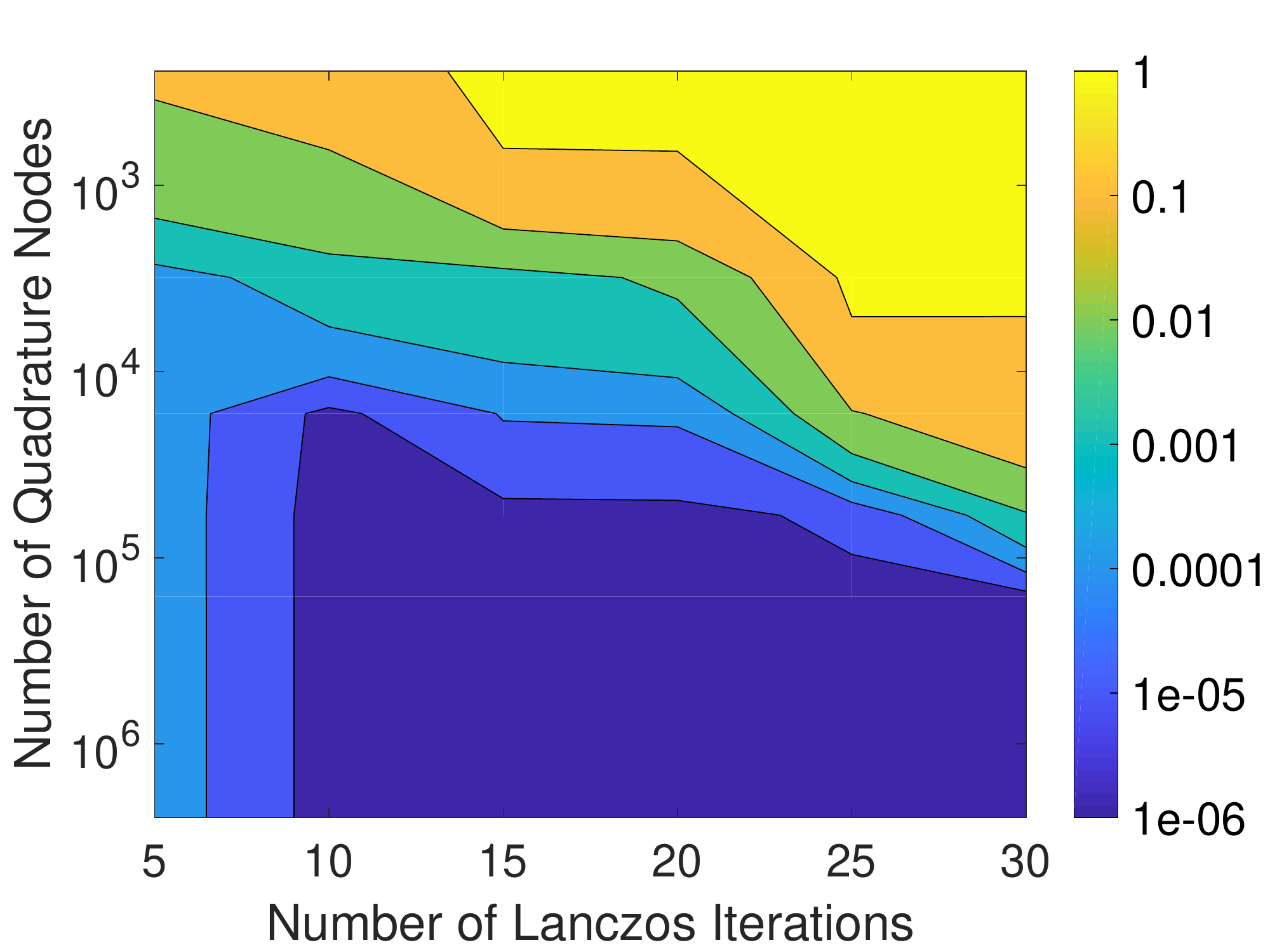}
}
\caption{Convergence study for Algorithm \ref{alg:LSAVE} on \eqref{eq:ex2}. Figure \ref{fig:LSAVE_quadpts_Ex2} depicts differences between subsequent matrices for increasing numbers of quadrature nodes with $k = 35$. Figure \ref{fig:LSAVE_lanciters_Ex2} contains matrix differences for increasing numbers of Lanczos iterations with 21 quadrature nodes per dimension ($N = 21^5 =$ 4,084,101). Figure \ref{fig:LSAVE_2d_converge_Ex2} shows the errors in the estimated $\CAVE$ matrix for various values of $N$ and $K$ relative to the ``true'' matrix (i.e., $N = 21^5 =$ 4,084,101, $k = 35$).}
\label{fig:LSAVE_relconv_Ex2}
\end{figure}

We next perform the same convergence study on the LSAVE algorithm for \eqref{eq:ex2}. Figures \ref{fig:LSAVE_quadpts_Ex2} and \ref{fig:LSAVE_lanciters_Ex2} contain the matrix difference study for increasing $N$ ($k = 35$ fixed) and increasing $k$ ($N = 21^5 =$ 4,084,101 fixed), respectively. Subsequent matrix differences decay as $N$ and $k$ increase, suggesting convergence of the Lanczos-Stieltjes approximation $\hCAVE$. Therefore, we use $N = 21^5 =$ 4,084,101 and $k = 35$ in Algorithm \ref{alg:LSAVE} to compute the ``true'' value of $\CAVE$. Figure \ref{fig:LSAVE_2d_converge_Ex2} contains the matrix errors from LSAVE for various values of $N$ and $k$ relative to this ``true'' $\CAVE$ matrix. The results of this study and their interpretations are similar to those of the LSIR study above.

\begin{figure}[!ht]
\centering
\subfloat[]{
\label{fig:slicestudy_SIR_Ex2}
\includegraphics[width=0.45\textwidth]{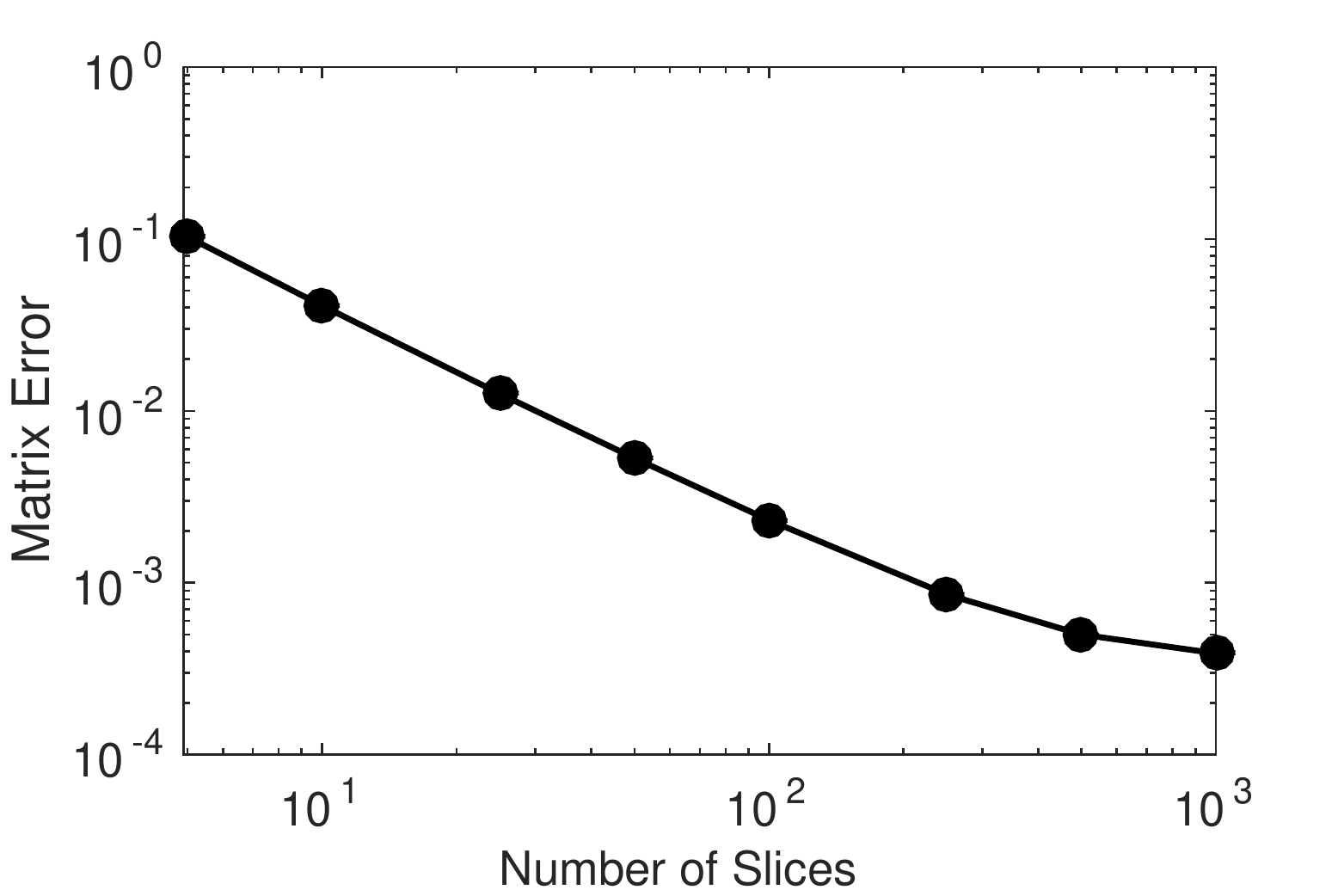}
}
\hfil
\subfloat[]{
\label{fig:slicestudy_SAVE_Ex2}
\includegraphics[width=0.45\textwidth]{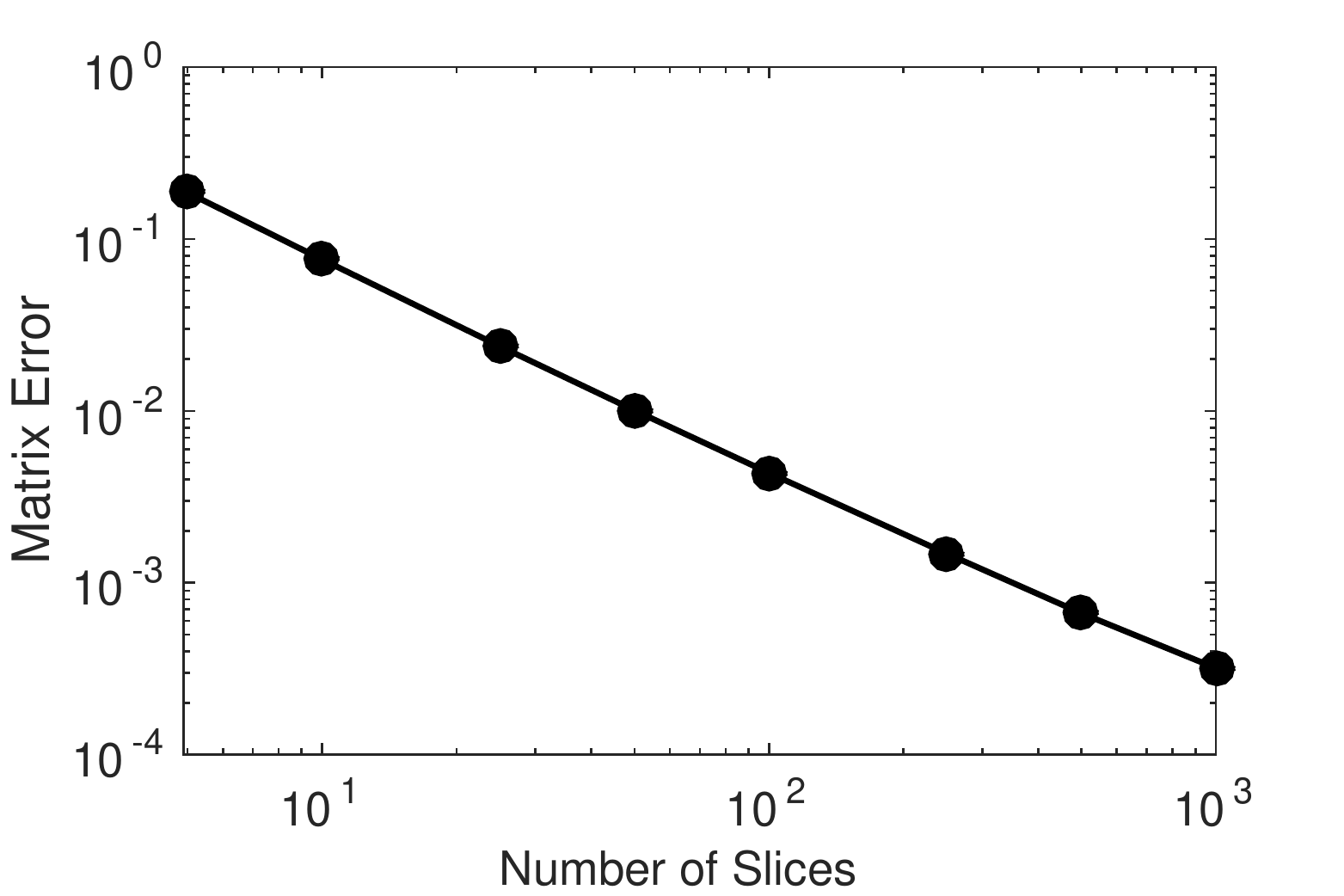}
}
\caption{A comparison of the Riemann sum approximation and the Lanczos-Stieltjes approximation of $\CIR$ (Figure \ref{fig:slicestudy_SIR_Ex2}) and $\CAVE$ (Figure \ref{fig:slicestudy_SAVE_Ex2}) for increasing number of Riemann sums (or slices) for \eqref{eq:ex2}. For the Lanczos-Stieltjes approximations, we used $N = 21^5 =$ 4,084,101 Gauss-Christoffel quadrature nodes and $k = 35$ Lanczos iterations. For the slice-based approximations, we used $N = 10^8$ Monte Carlo samples.}
\label{fig:slicestudy_Ex2}
\end{figure}

Next, we examine how the approximated SIR and SAVE matrices from Algorithms \ref{alg:SIR} and \ref{alg:SAVE}, respectively, compare to the LSIR and LSAVE approximations. Recall $\hCSIR$ and $\hCSAVE$ contain two levels of approximation---one due to the number of samples $N$ and one due to the number of Riemann sums (or slices) $R$ over the output space. We first focus on convergence in terms of Riemann sums. Figure \ref{fig:slicestudy_Ex2} compares the approximated SIR and SAVE matrices for increasing $R$ to their Lanczos-Stieltjes counterparts. For the Lanczos-Stieltjes approximations---$\hCIR$ and $\hCAVE$---we use $N = 21^5 =$ 4,084,101 quadrature nodes and $k = 35$ Lanczos iterations as this was shown to produce sufficiently converged matrices in the previous study. For the SIR and SAVE algorithms, we use $N = 10^8$ samples randomly drawn according to $\rho$. Additionally, we define the slices adaptively based on the sampling to balance the number of samples in each slice as discussed in Section \ref{sec:IR_methods}. We notice that the sliced approximations converge to their Lanczos-Stieltjes counterparts at a rate $R^{-1}$ as expected for Riemann sums~\citep[Ch. 2]{Davis84}.

Finally, we compare the slice-based algorithms to their Lanczos-Stieltjes counterparts in terms of sampling. The SIR and SAVE algorithms can be interpreted as computing Monte Carlo approximations of the conditional expectation and covariance in the context of computer experiments. However, the LSIR and LSAVE algorithms use tensor product Gauss-Christoffel quadrature rules with respect to $\rho$. To enable direct comparison, we perform Algorithms \ref{alg:LSIR} and \ref{alg:LSAVE} using $\vx_i$ randomly sampled according to $\rho$ with weights $\nu_i = 1/N$ for $i = 0,\dots,N-1$. Note that this Monte Carlo integration rule can be interpreted as introducing a discrete norm and inner product as discussed in Section \ref{subsec:lanczos}. We perform this comparison for various values of $R$ (the number of Riemann sums or slices) and  $k$ (the number of Lanczos iterations) for each of the methods. Based on the previous studies, we again use the Lanczos-Stieltjes algorithms with $N = 21^5 =$ 4,084,101 quadrature nodes and $k = 35$ Lanczos iterations as the ``true'' values of $\CIR$ and $\CAVE$ for computation of the relative matrix errors.

\begin{figure}[!ht]
\centering
\subfloat[]{
\label{fig:SIR_samples_MC_Ex2}
\includegraphics[width=0.45\textwidth]{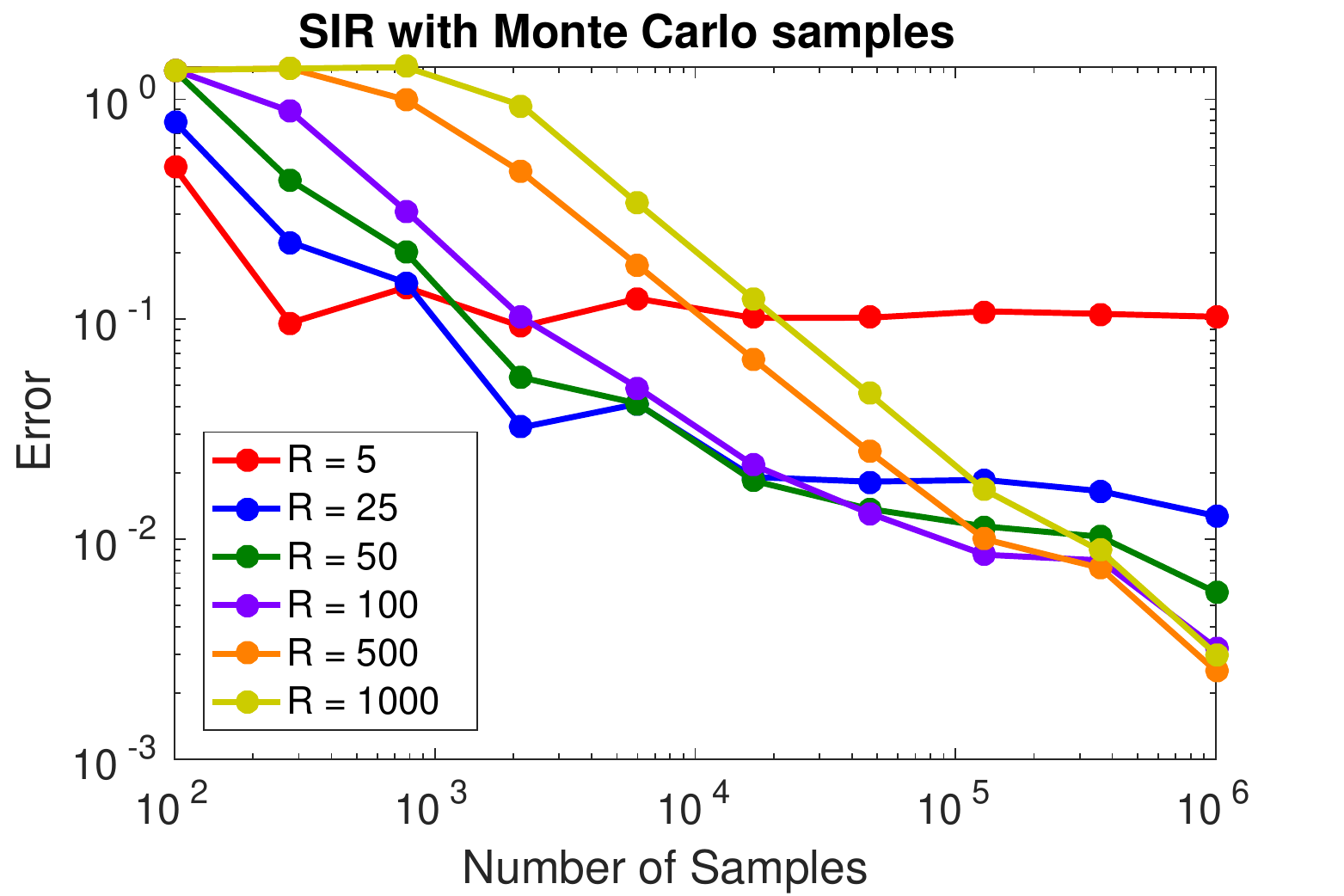}
}
\hfil
\subfloat[]{
\label{fig:LSIR_samples_MC_Ex2}
\includegraphics[width=0.45\textwidth]{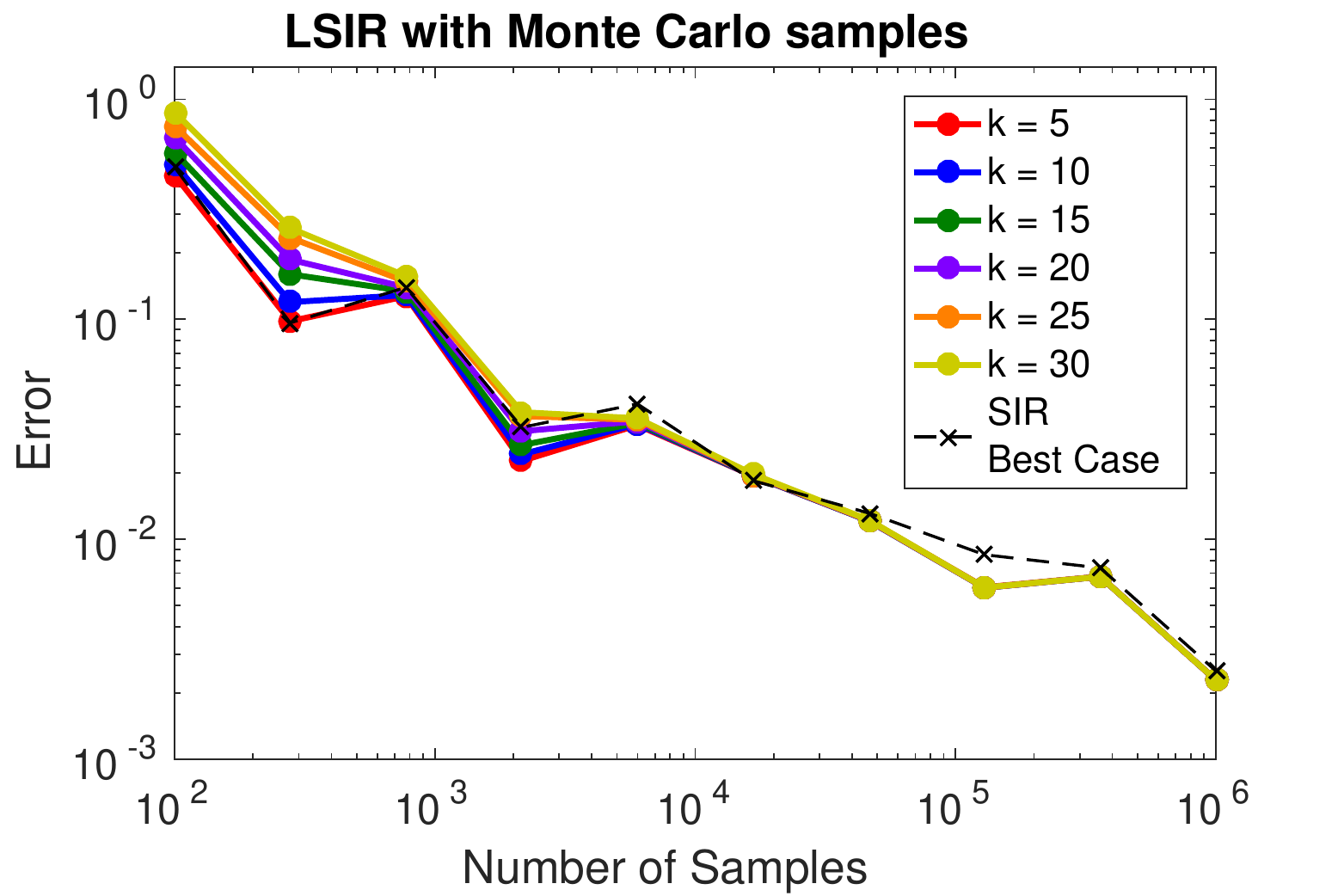}
}
\caption{A comparison of the SIR and LSIR algorithms for \eqref{eq:ex2} using Monte Carlo samples. Figure \ref{fig:SIR_samples_MC_Ex2} shows the relative matrix errors as a function of the number of samples for various values of $R$ (the number of slices). Figure \ref{fig:LSIR_samples_MC_Ex2} contains the relative matrix errors as a function of the number of samples for various values of $k$ (the number of Lanczos iterations). This plot also contains the best case result from the SIR algorithm for reference.}
\label{fig:LSIR_samples_Ex2}
\end{figure}

Figure \ref{fig:LSIR_samples_Ex2} contains the results of this study for the SIR and LSIR algorithms on \eqref{eq:ex2}. The first plot contains relative errors in $\hCSIR$ for increasing numbers of Monte Carlo samples at various levels of slicing. The second plot shows the relative errors for the LSIR algorithm using the same Monte Carlo samples as in the SIR case for various numbers of Lanczos iterations. We notice less variance in the LSIR plot due to $k$ than in the SIR plot due to $R$. Additionally, the slicing in \ref{fig:SIR_samples_MC_Ex2} is performed adaptively to balance the numbers of samples within each slice. The issue of choosing how to slice up the output space does not exist in the Lanczos-Stieltjes approach as the method automatically chooses the best integration rule over $\sF$. Figure \ref{fig:LSIR_samples_MC_Ex2} also includes the best case error from the SIR plot. That is, it plots the minimum error among all of the tested values of $R$ for each value of $N$. This shows the LSIR algorithm performing approximately as well as the best case in SIR, regardless of the value of $k$ chosen. This suggests that LSIR not only enables the use of high-order quadrature when the dimension is sufficiently small, but it can also be used with Monte Carlo sampling and perform as well as SIR. 

\begin{figure}[!ht]
\centering
\subfloat[]{
\label{fig:SAVE_samples_MC_Ex2}
\includegraphics[width=0.45\textwidth]{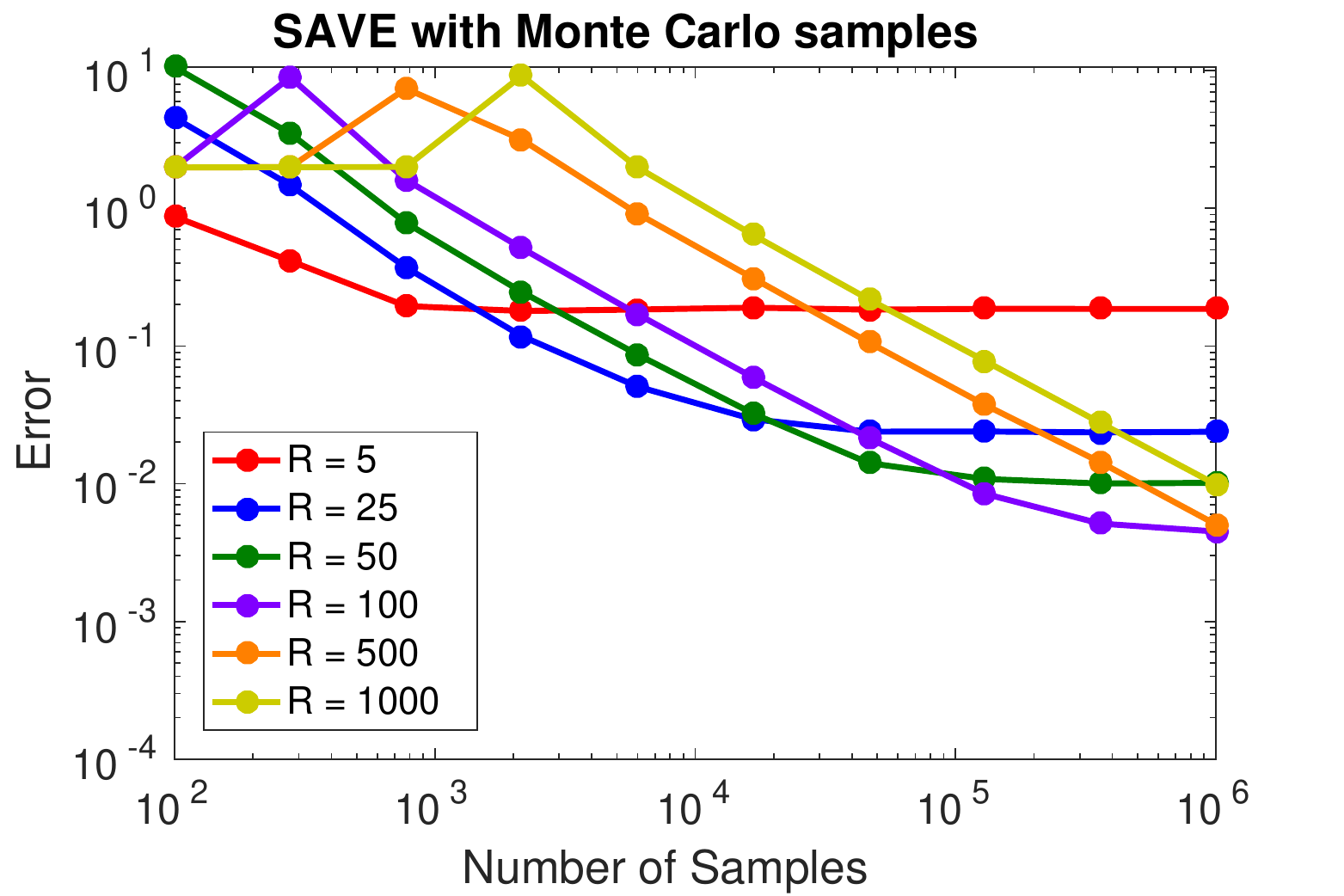}
}
\hfil
\subfloat[]{
\label{fig:LSAVE_samples_MC_Ex2}
\includegraphics[width=0.45\textwidth]{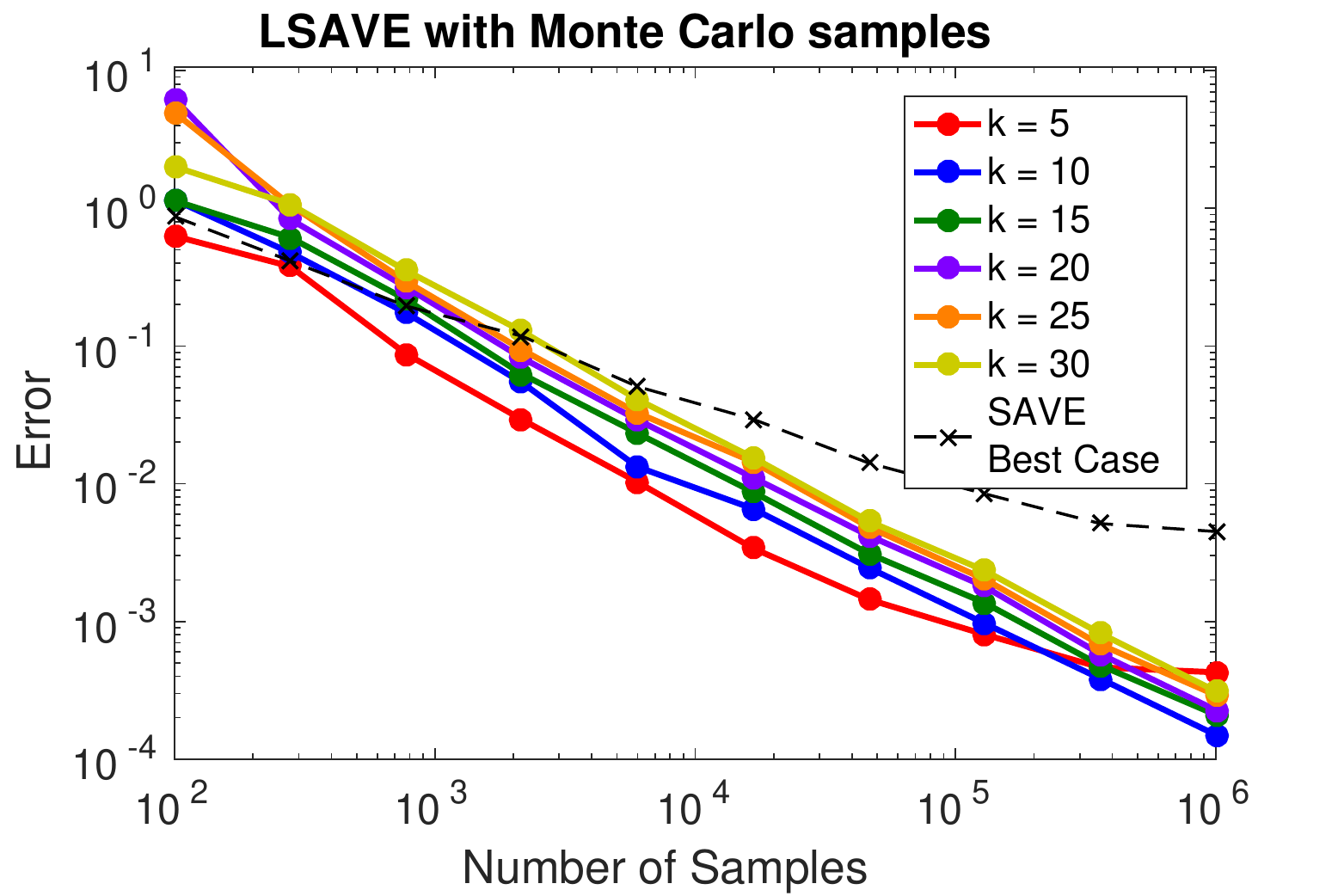}
}
\caption{A comparison of the SAVE and LSAVE algorithms for \eqref{eq:ex2} using Monte Carlo samples. Figure \ref{fig:SAVE_samples_MC_Ex2} shows the relative matrix errors as a function of the number of samples for various values of $R$ (the number of slices). Figure \ref{fig:LSAVE_samples_MC_Ex2} contains the relative matrix errors as a function of the number of samples for various values of $k$ (the number of Lanczos iterations). This plot also contains the best case result from the SAVE algorithm for reference.}
\label{fig:LSAVE_samples_Ex2}
\end{figure}

Figure \ref{fig:LSAVE_samples_Ex2} contains the results of the same study as above but applied to the SAVE and LSAVE algorithms. The discussion of the results mirrors that of Figure \ref{fig:LSIR_samples_Ex2}. One notable difference is that the LSAVE algorithm outperforms the best case SAVE results for large $N$. This again suggests that the Lanczos-Stieltjes approach can be used even when the dimension of the problem prohibits the use of tensor product quadrature rules.


\subsection{Example 3: OTL circuit function}
\label{subsec:OTL_circuit}

In this section, we consider the physically-motivated OTL circuit function. This function models the midpoint voltage ($V_m$) from a transformerless push-pull circuit (see~\citet{Surjanovic15} and~\citet{BenAri07}). It has the form
\begin{equation} \label{eq:ex3}
V_m
\;=\;
\frac{\left( \left( 12 R_{b2} (R_{b1} + R_{b2})^{-1} + 0.74 \right) \beta \left( R_{c2} + 9 \right) + 11.35 R_{f} \right) R_{c1} + 0.74 R_{f} \beta \left( R_{c2} + 9 \right)}{\left( \beta \left( R_{c2} + 9 \right) + R_{f} \right) R_{c1}} ,
\end{equation}
where the 6 physical inputs are described in Table \ref{tab:OTL_inputs}. The table also contains the ranges of the inputs. We assume a uniform density applied over the the ranges given in the table.
\begin{table}[!ht]
\centering
\caption{The 6 physical inputs for the OTL circuit function with their ranges. We assume a uniform density over these ranges.}
\label{tab:OTL_inputs}
\begin{tabular}{ccc}
Input & Description (Units) & Range \\ \hline
$R_{b1}$ & resistance b1 (K-Ohms) & $[50, 150]$ \\
$R_{b2}$ & resistance b2 (K-Ohms) & $[25, 70]$ \\
$R_{f}$ & resistance f (K-Ohms)& $[0.5, 3]$ \\
$R_{c1}$ & resistance c1 (K-Ohms)& $[1.2, 2.5]$ \\
$R_{c2}$ & resistance c2 (K-Ohms)& $[0.25, 1.2]$ \\
$\beta$ & current gain (Amperes) & $[50, 300]$
\end{tabular}
\end{table}

\begin{figure}[!ht]
\centering
\subfloat[]{
\label{fig:LSIR_quadpts_Ex3}
\includegraphics[width=0.4\textwidth]{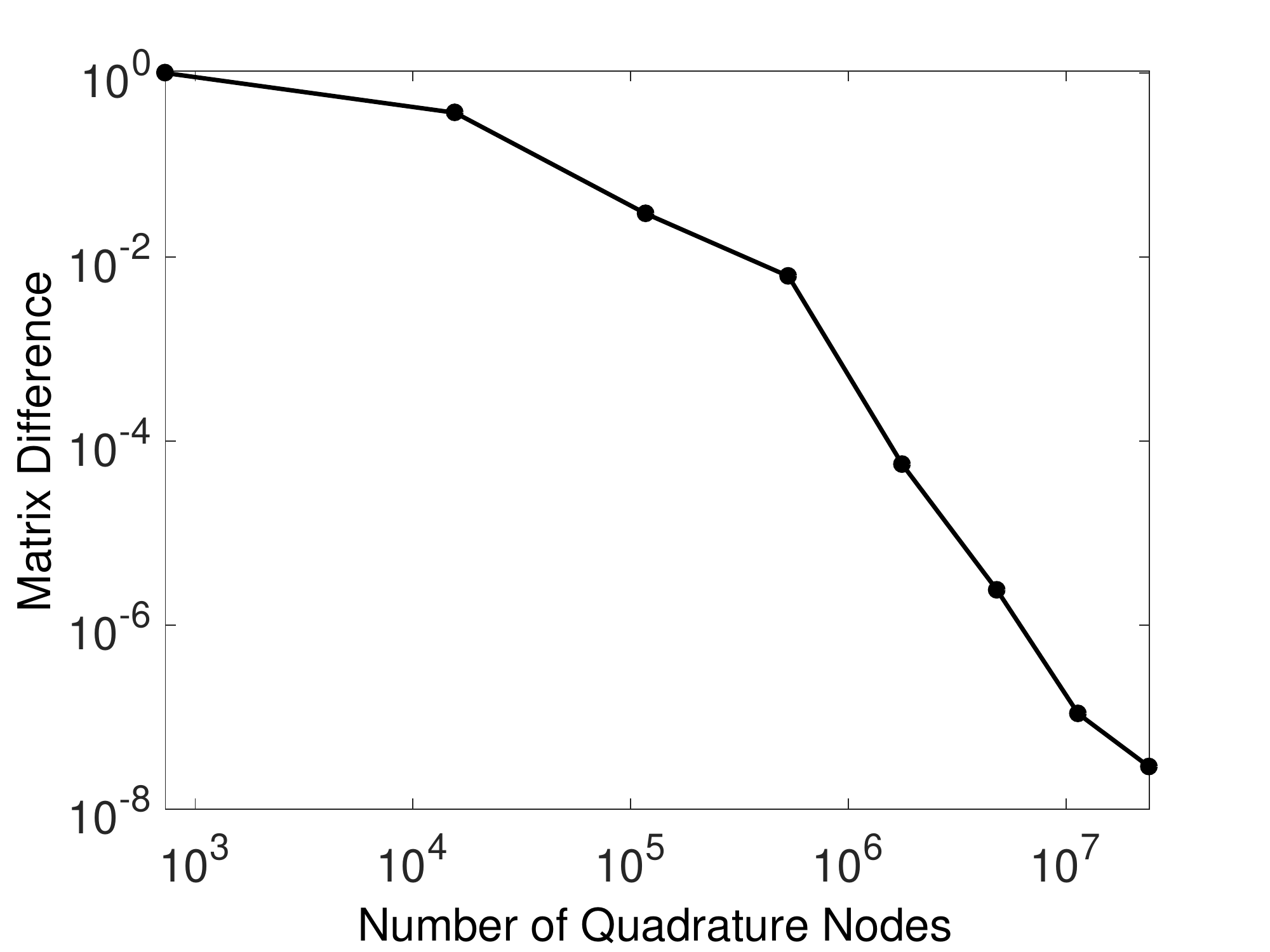}
}
\subfloat[]{
\label{fig:LSIR_lanciters_Ex3}
\includegraphics[width=0.4\textwidth]{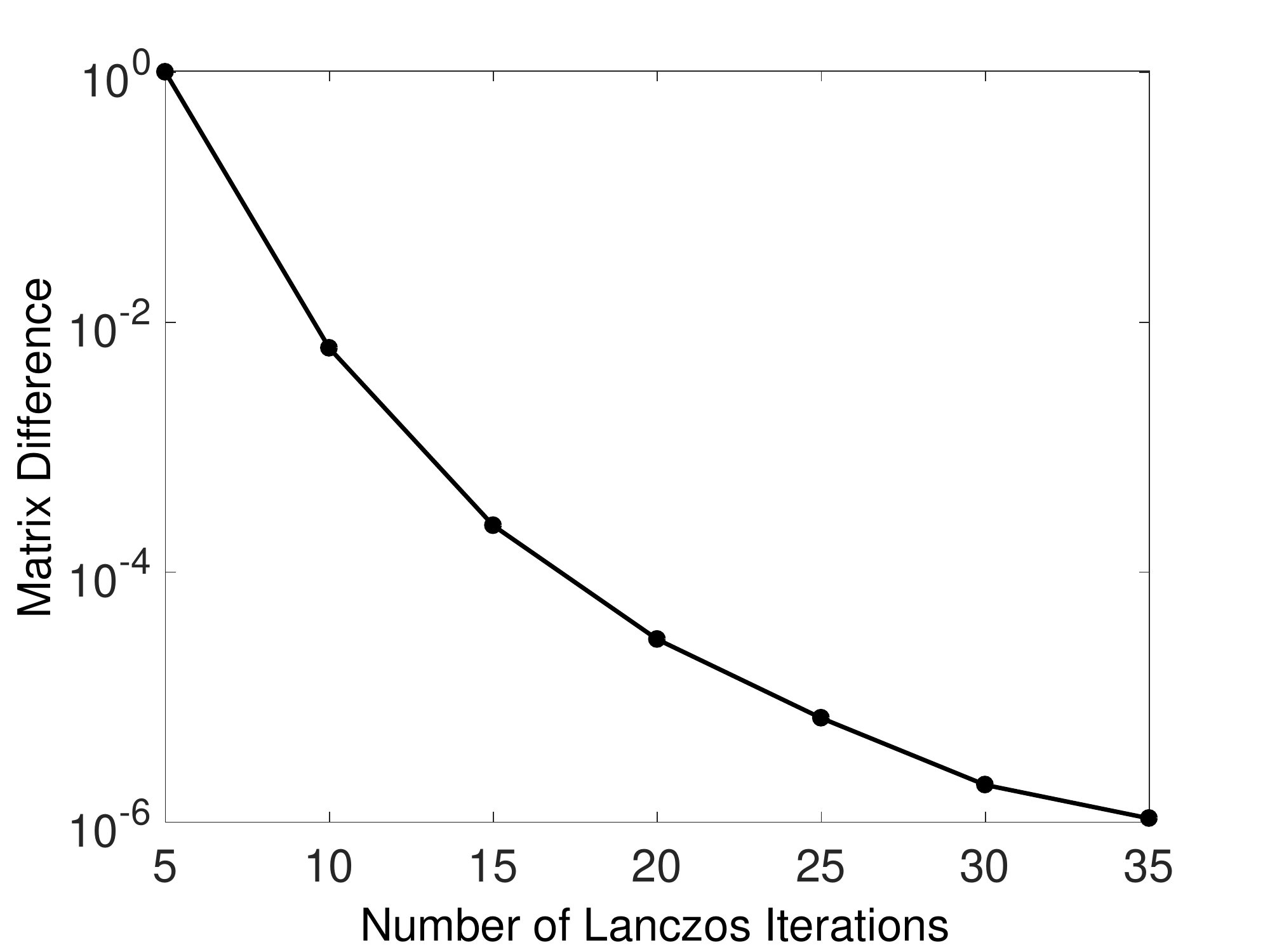}
} \\
\subfloat[]{
\label{fig:LSIR_2d_converge_Ex3}
\includegraphics[width=0.4\textwidth]{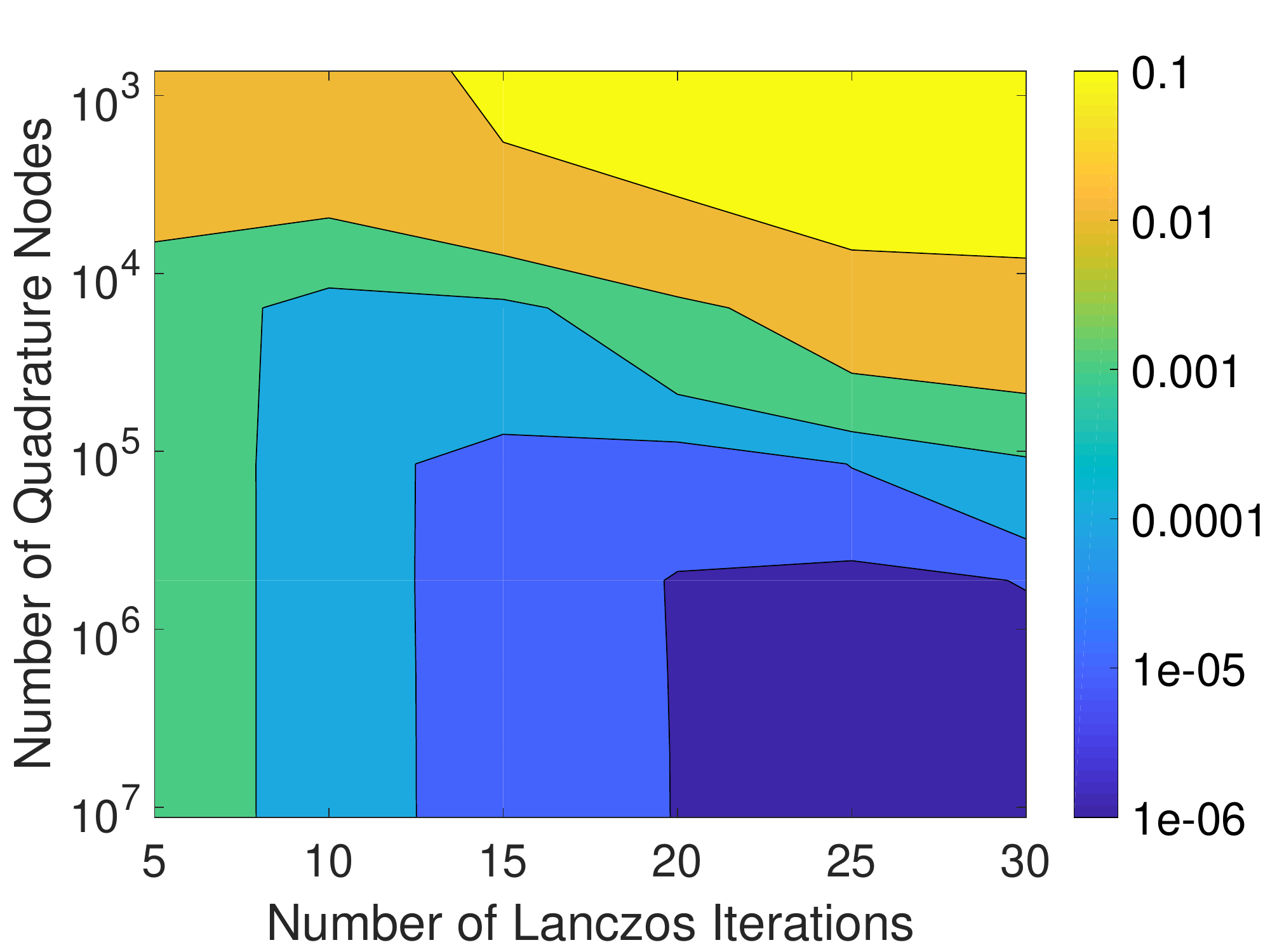}
}
\caption{Convergence study for Algorithm \ref{alg:LSIR} on \eqref{eq:ex3}. Figure \ref{fig:LSIR_quadpts_Ex3} depicts differences between subsequent matrices for increasing numbers of quadrature nodes with $k = 35$. Figure \ref{fig:LSIR_lanciters_Ex3} contains matrix differences for increasing numbers of Lanczos iterations with 21 quadrature nodes per dimension ($N = 17^6 =$ 24,137,569). Figure \ref{fig:LSIR_2d_converge_Ex3} shows the errors in the estimated $\CIR$ matrix for various values of $N$ and $K$ relative to the ``true'' matrix (i.e., $N = 17^6 =$ 24,137,569, $k = 35$).}
\label{fig:LSIR_relconv_Ex3}
\end{figure}

\begin{figure}[!ht]
\centering
\subfloat[]{
\label{fig:LSAVE_quadpts_Ex3}
\includegraphics[width=0.4\textwidth]{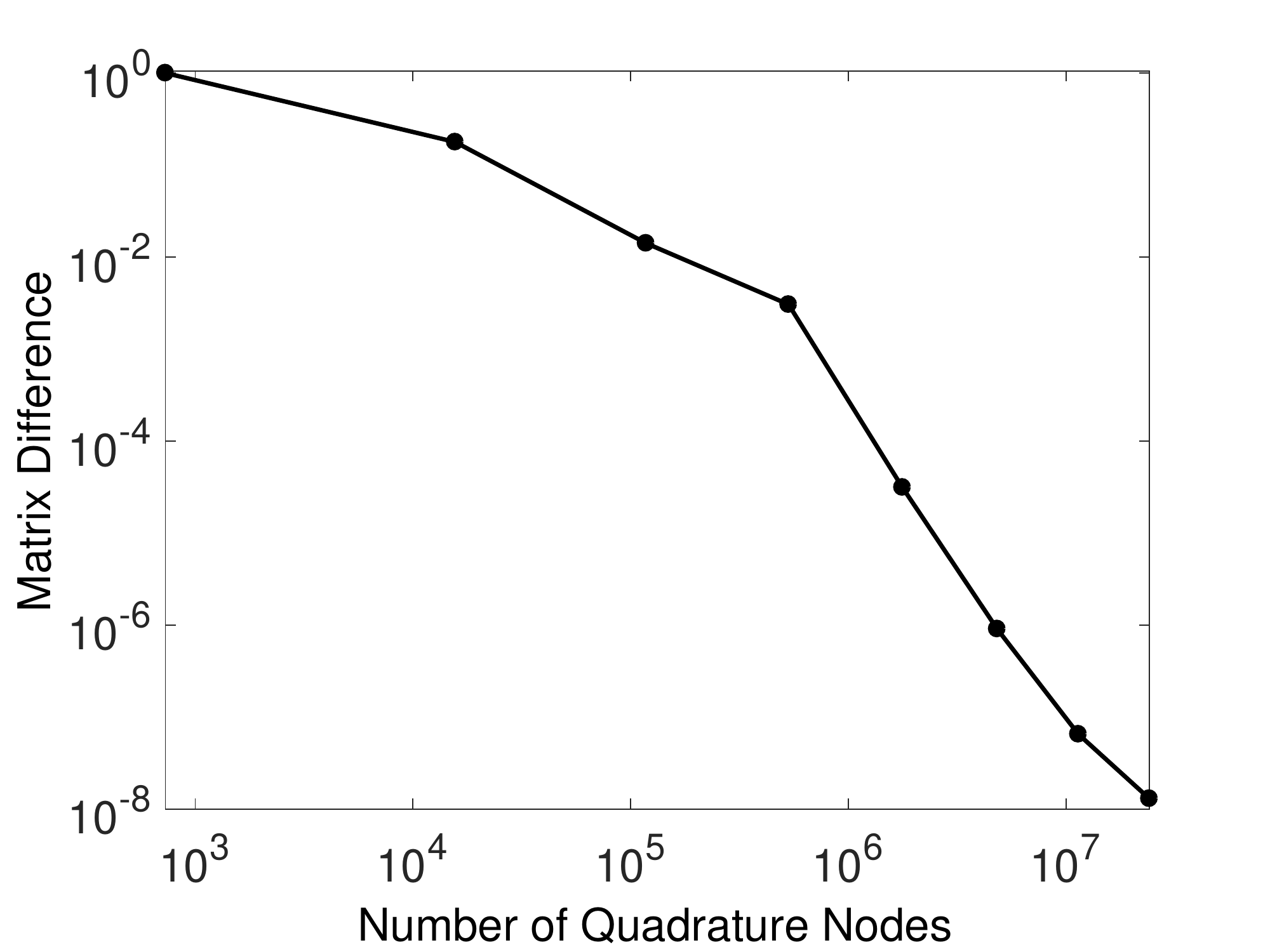}
}
\subfloat[]{
\label{fig:LSAVE_lanciters_Ex3}
\includegraphics[width=0.4\textwidth]{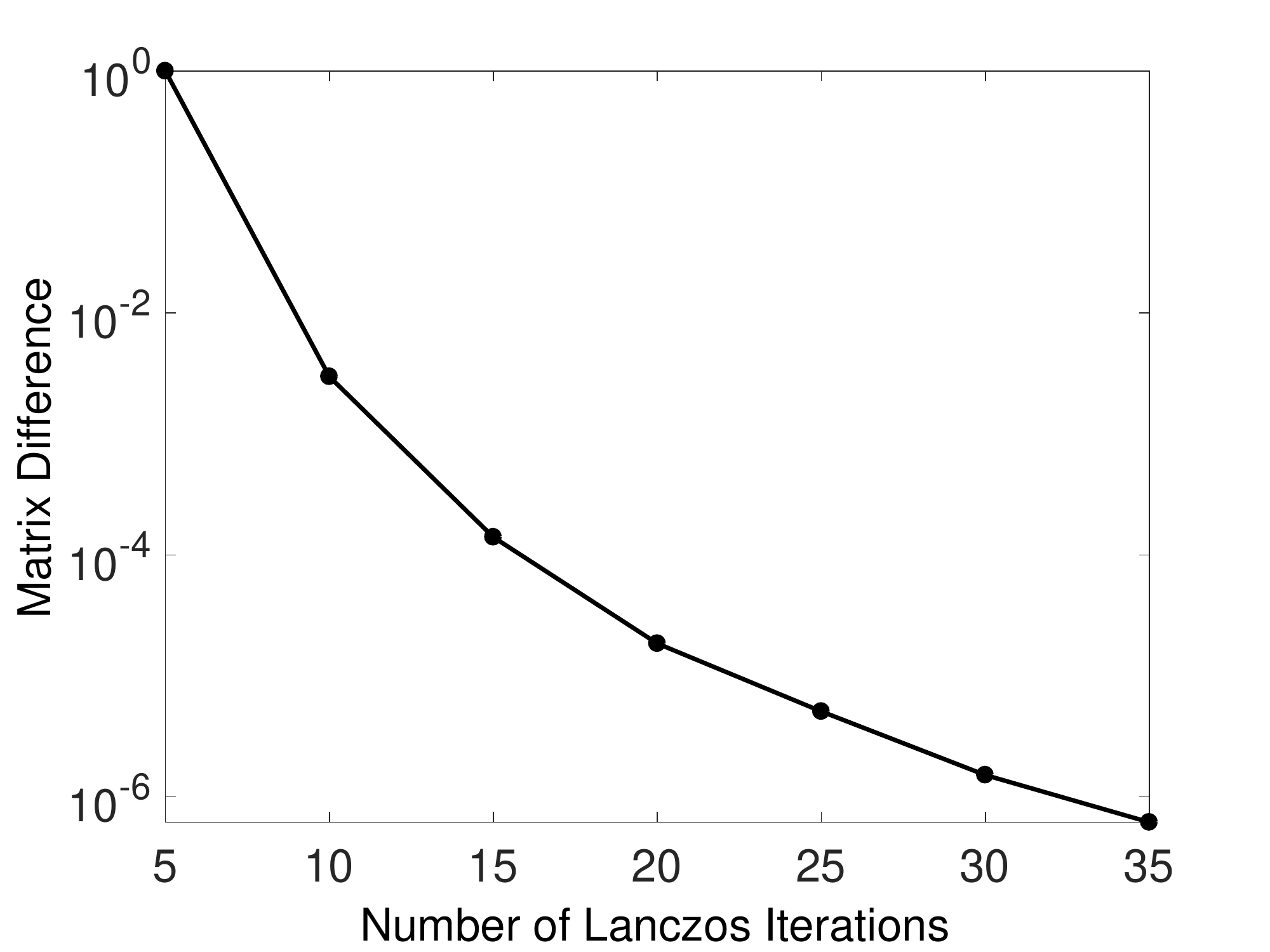}
} \\
\subfloat[]{
\label{fig:LSAVE_2d_converge_Ex3}
\includegraphics[width=0.4\textwidth]{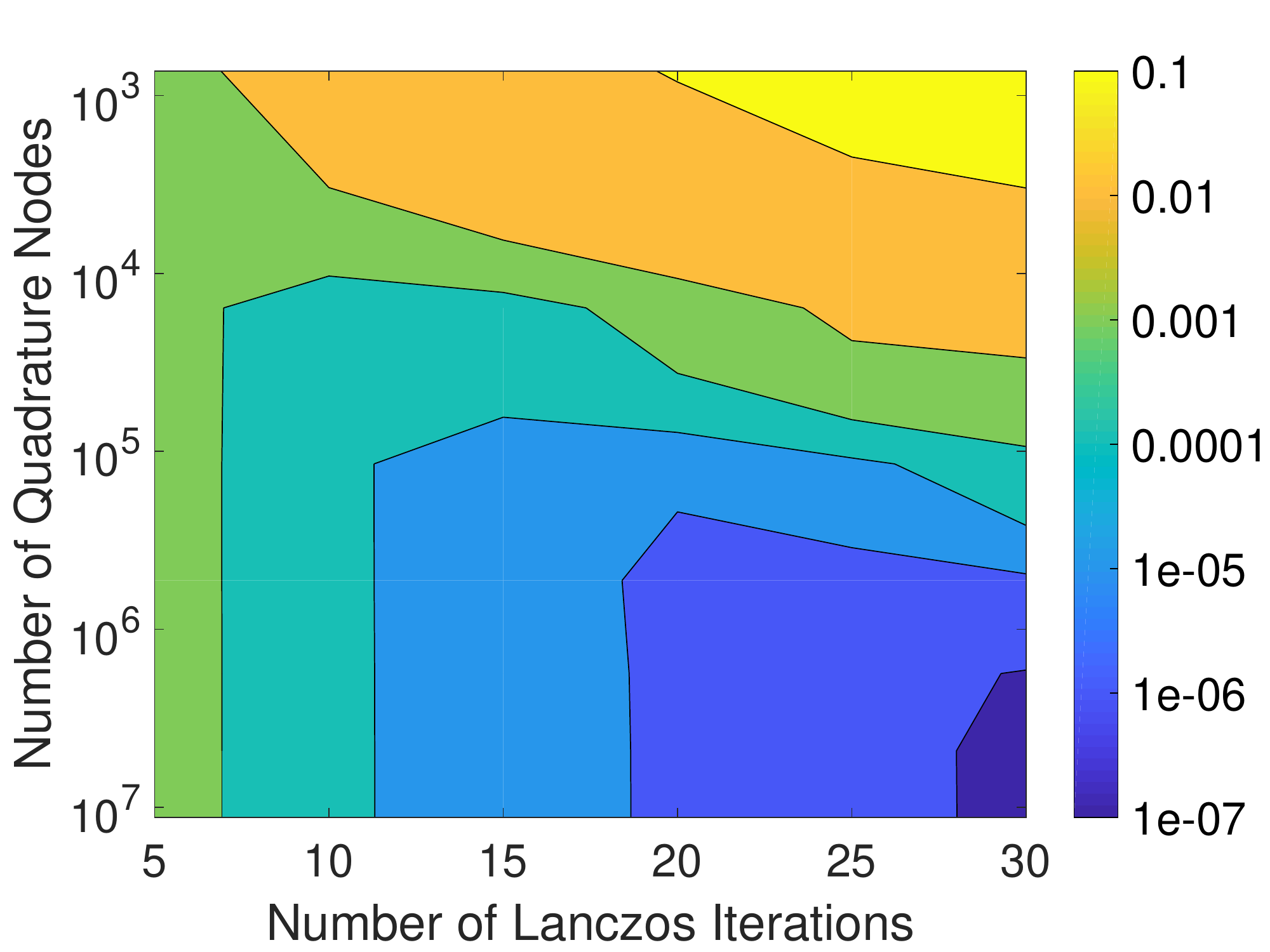}
}
\caption{Convergence study for Algorithm \ref{alg:LSAVE} on \eqref{eq:ex3}. Figure \ref{fig:LSAVE_quadpts_Ex3} depicts differences between subsequent matrices for increasing numbers of quadrature nodes with $k = 35$. Figure \ref{fig:LSAVE_lanciters_Ex3} contains matrix differences for increasing numbers of Lanczos iterations with 21 quadrature nodes per dimension ($N = 17^6 =$ 24,137,569). Figure \ref{fig:LSAVE_2d_converge_Ex3} shows the errors in the estimated $\CAVE$ matrix for various values of $N$ and $K$ relative to the ``true'' matrix (i.e., $N = 17^6 =$ 24,137,569, $k = 35$).}
\label{fig:LSAVE_relconv_Ex3}
\end{figure}

We perform the same convergence study from Section \ref{subsec:1d_ridge} for the LSIR and LSAVE algorithms with various numbers of Gauss-Christoffel quadrature nodes $N$ and Lanczos iterations $k$ applied to the OTL circuit function. Figure \ref{fig:LSIR_relconv_Ex3} contains the results of this study for LSIR. Figure \ref{fig:LSIR_quadpts_Ex3} shows the subsequent matrix differences for increasing $N$ with $k = 35$, and Figure \ref{fig:LSIR_lanciters_Ex3} shows the subsequent matrix differences for increasing $k$ with 17 quadrature nodes per dimension (i.e., $N = 17^6 =$ 24,137,569). Both of these plots suggest that the Lanczos-Stieltjes estimate of $\CIR$ is converging in terms of $N$ and $k$ for \eqref{eq:ex3}. Thus, we define the ``true'' $\CIR$ matrix to be the result of Algorithm \ref{alg:LSIR} applied to \eqref{eq:ex3} with $N = 17^6 =$ 24,137,569 and $k = 35$. Figure \ref{fig:LSIR_2d_converge_Ex3} contains matrix errors relative to this ``true'' $\CIR$ for various values of $N$ and $k$. We see similar behavior as in Section \ref{subsec:1d_ridge}. Figure \ref{fig:LSAVE_relconv_Ex3} contains the results of this study for LSAVE.

\begin{figure}[!ht]
\centering
\subfloat[]{
\label{fig:slicestudy_SIR_Ex3}
\includegraphics[width=0.45\textwidth]{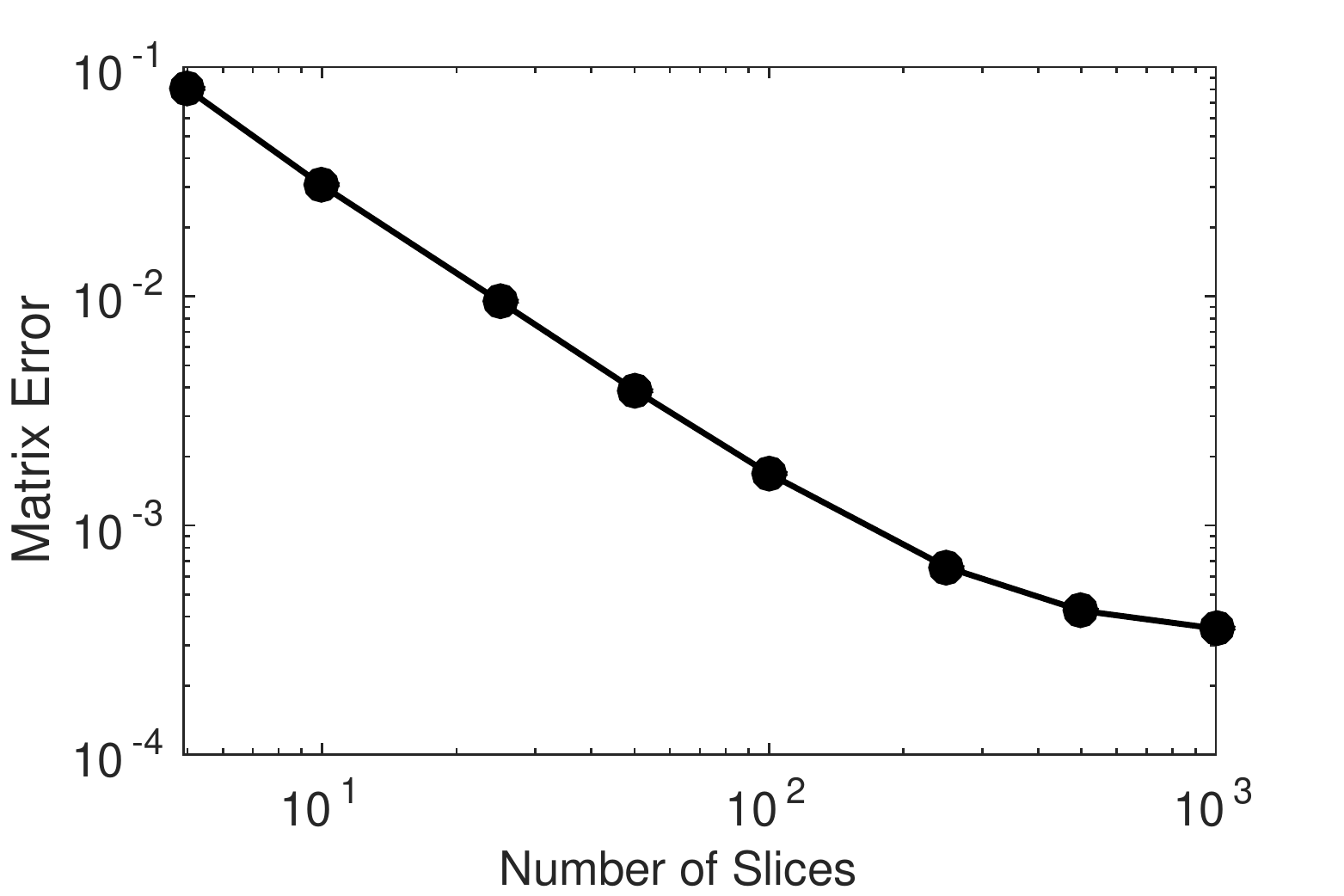}
}
\hfil
\subfloat[]{
\label{fig:slicestudy_SAVE_Ex3}
\includegraphics[width=0.45\textwidth]{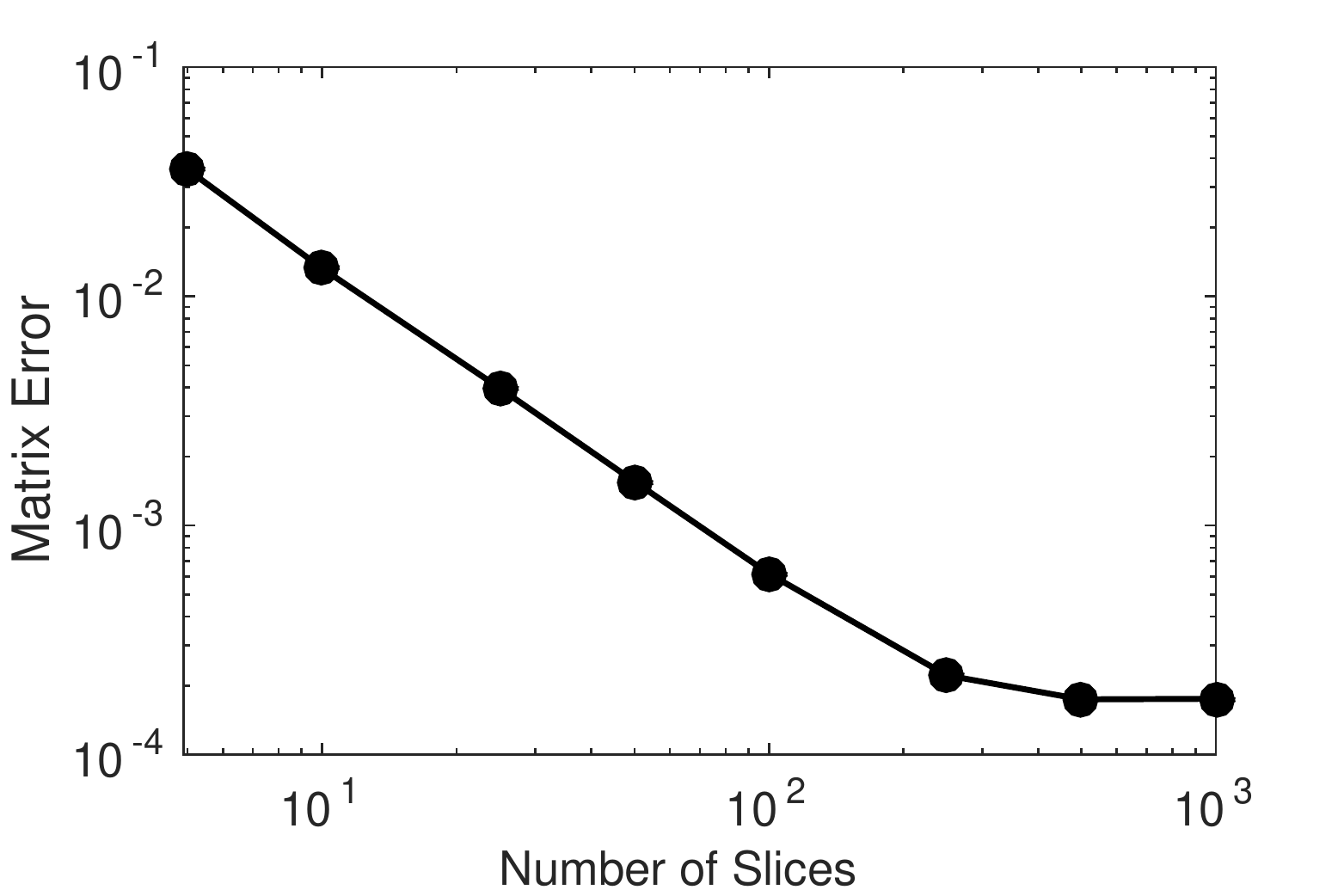}
}
\caption{A comparison of the Riemann sum approximation and the Lanczos-Stieltjes approximation of $\CIR$ (Figure \ref{fig:slicestudy_SIR_Ex3}) and $\CAVE$ (Figure \ref{fig:slicestudy_SAVE_Ex3}) for increasing number of Riemann sums (or slices) for \eqref{eq:ex3}. For the Lanczos-Stieltjes approximations, we used $N = 21^5 =$ 4,084,101 Gauss-Christoffel quadrature nodes on $\sX$ and $k = 35$ Lanczos iterations. For the slice-based approximations, we used $N = 10^8$ Monte Carlo samples.}
\label{fig:slicestudy_Ex3}
\end{figure}

Lastly, we compare the Lanczos-Stieltjes algorithms to their sliced counterparts for the OTL circuit function. Figure \ref{fig:slicestudy_Ex3} plots matrix errors from SIR and SAVE for increasing number of slices $R$ and $N = 10^8$. The errors are computed relative to the Lanczos-Stieltjes matrices with $N = 17^6 =$ 24,137,569, $k = 35$. Again, the errors decay like $R^{-1}$ as expected for Riemann sum approximations.

\begin{figure}[!ht]
\centering
\subfloat[]{
\label{fig:SIR_samples_MC_Ex3}
\includegraphics[width=0.45\textwidth]{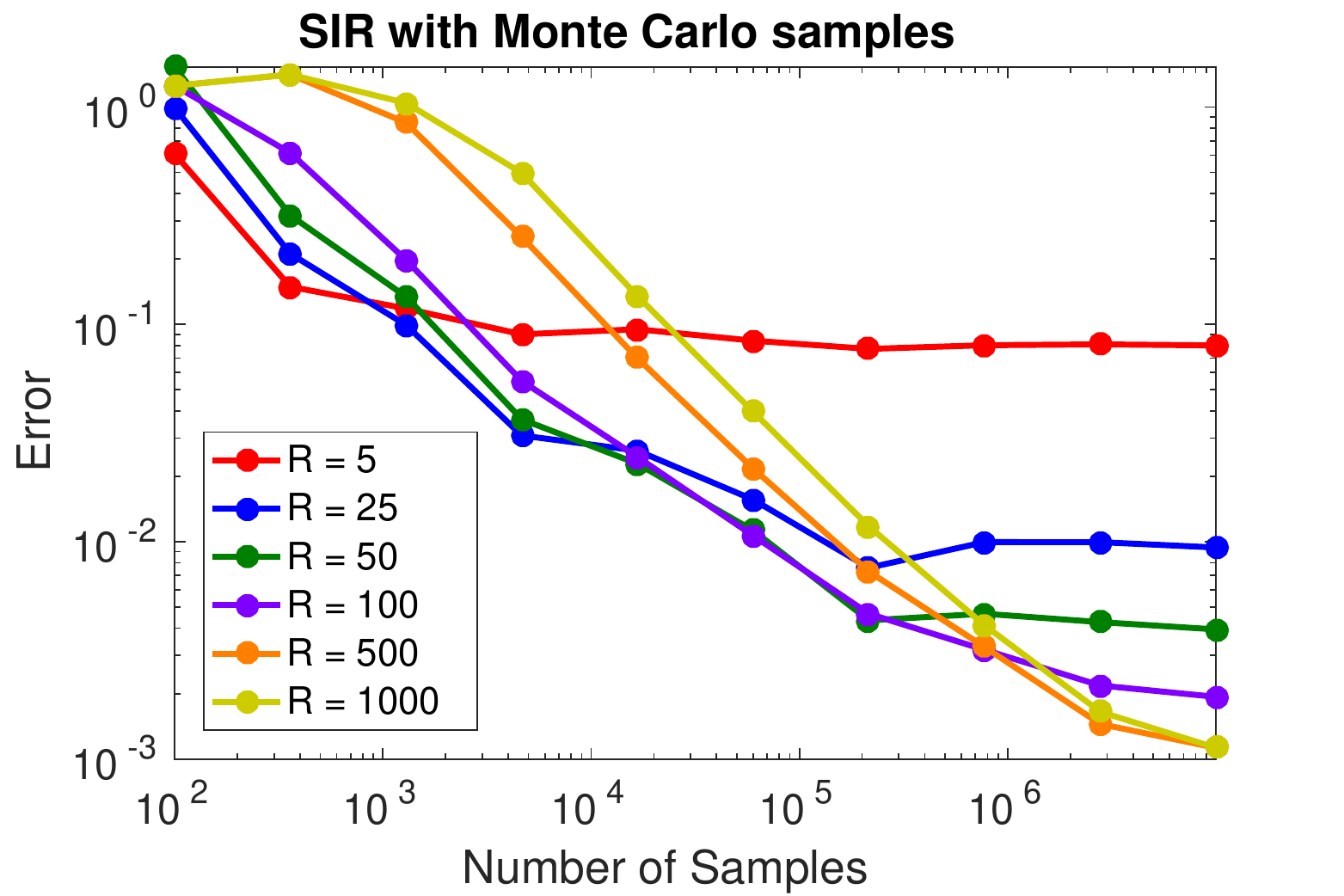}
}
\hfil
\subfloat[]{
\label{fig:LSIR_samples_MC_Ex3}
\includegraphics[width=0.45\textwidth]{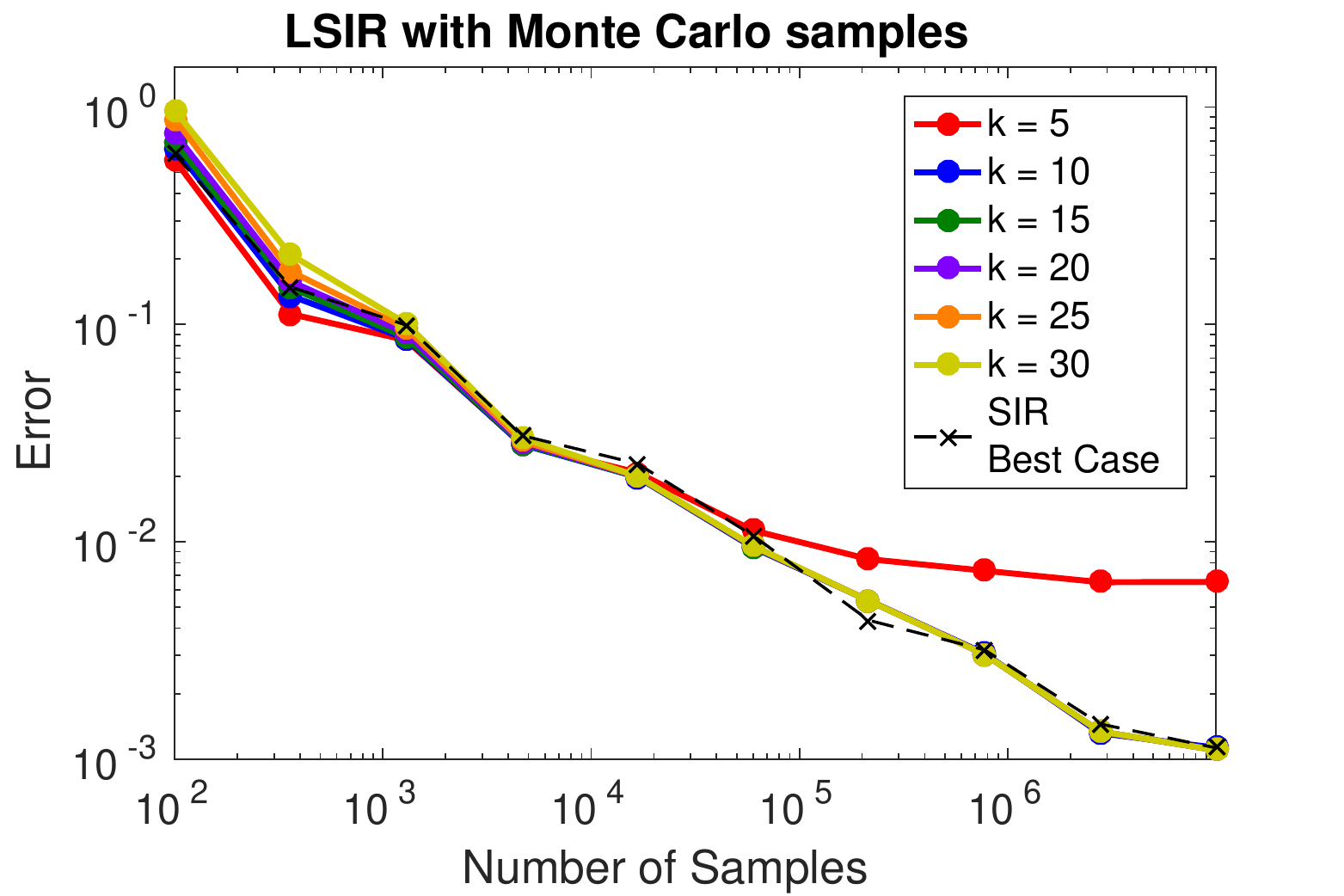}
} \\
\subfloat[]{
\label{fig:SAVE_samples_MC_Ex3}
\includegraphics[width=0.45\textwidth]{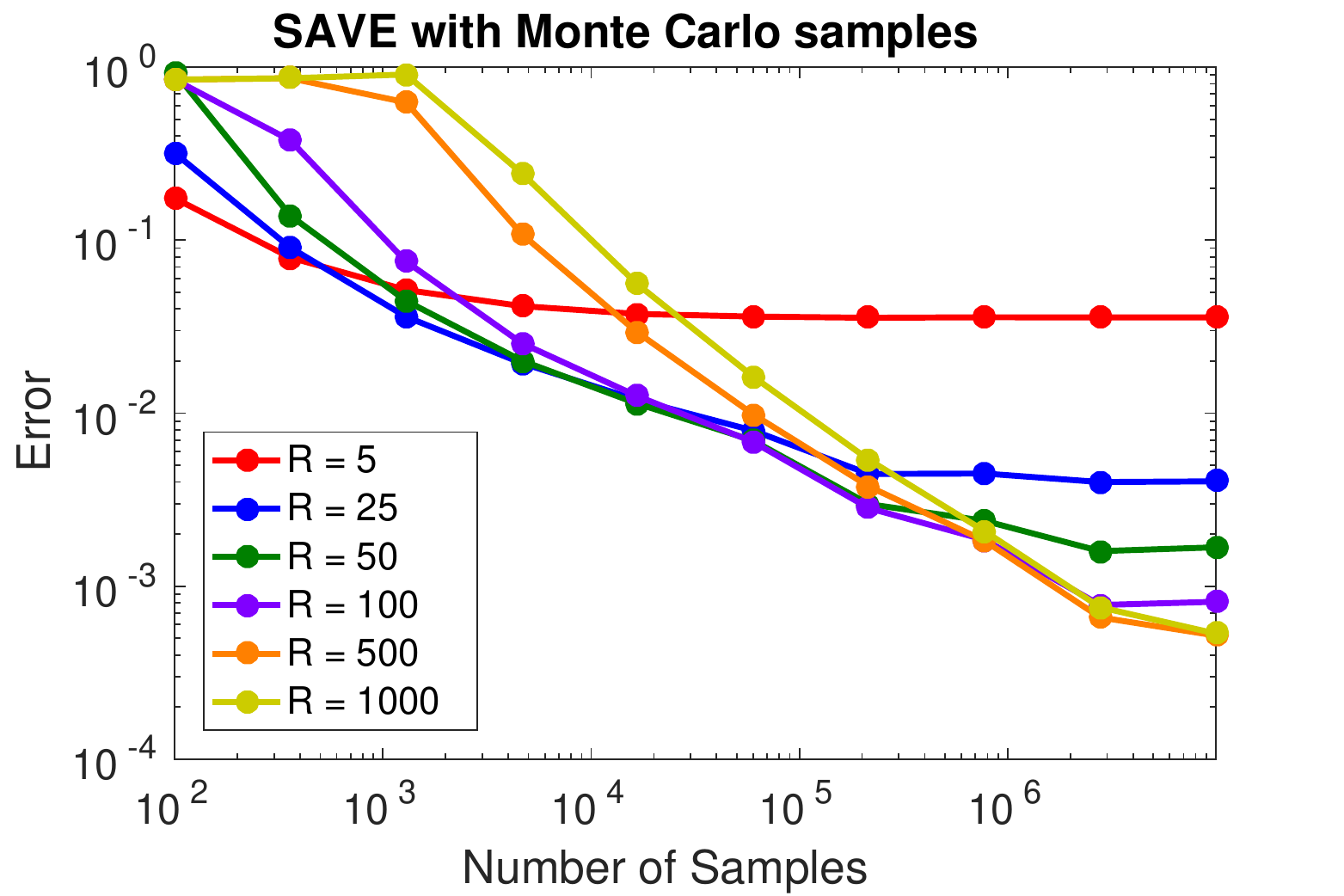}
}
\hfil
\subfloat[]{
\label{fig:LSAVE_samples_MC_Ex3}
\includegraphics[width=0.45\textwidth]{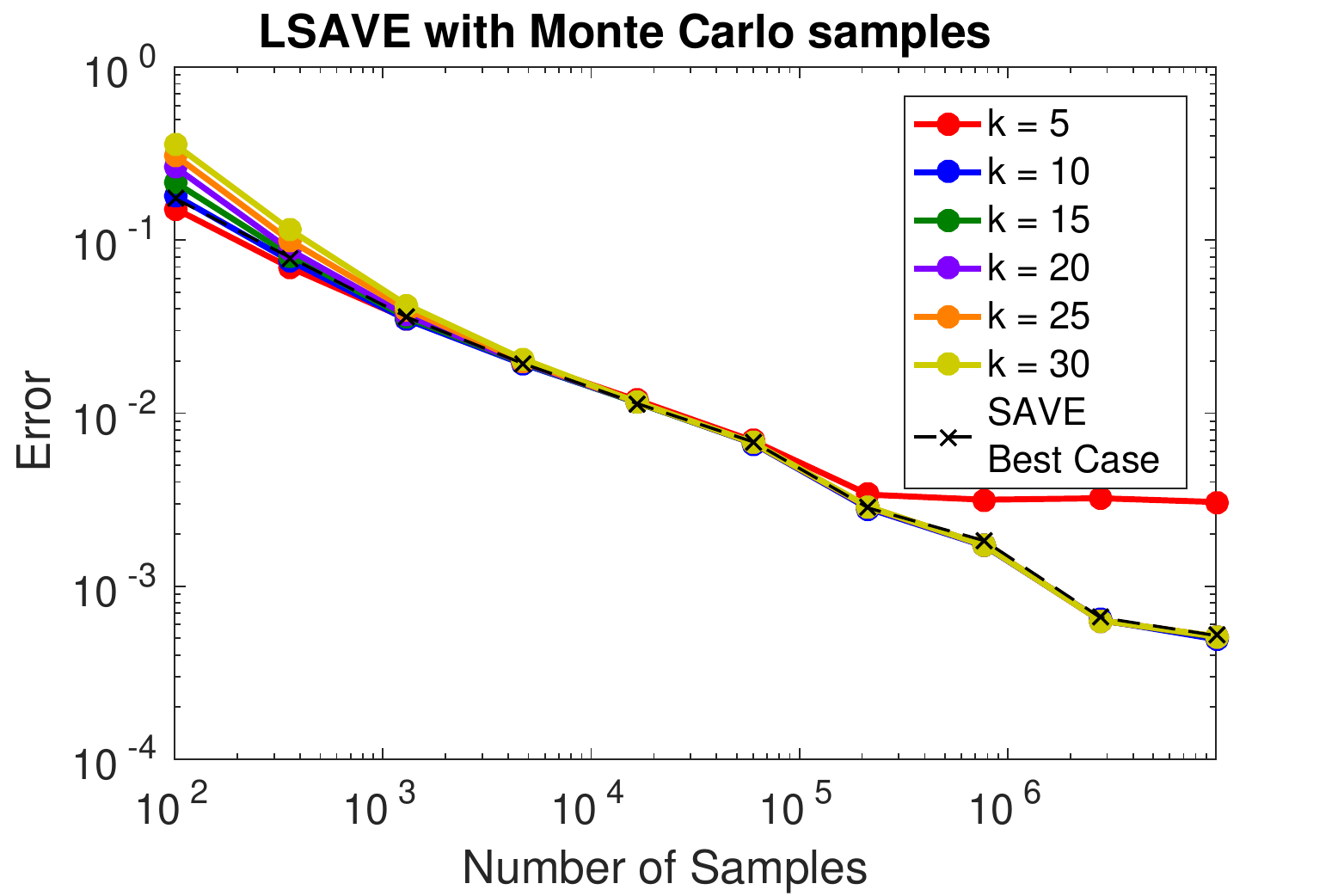}
}
\caption{A comparison of the SIR/SAVE and LSIR/LSAVE algorithms for \eqref{eq:ex3} using Monte Carlo samples. Figures \ref{fig:SIR_samples_MC_Ex3} and \ref{fig:SAVE_samples_MC_Ex3} show the relative matrix errors of the SIR and SAVE algorithms, respectively, as a function of the number of samples for various values of $R$ (the number of slices). Figures \ref{fig:LSIR_samples_MC_Ex3} and \ref{fig:LSAVE_samples_MC_Ex3} contain the relative matrix errors of the LSIR and LSAVE algorithms, respectively, as a function of the number of samples for various values of $k$ (the number of Lanczos iterations). These plots also contain the best case results from their slice-based counterparts for reference.}
\label{fig:samples_Ex3}
\end{figure}

Figure \ref{fig:samples_Ex3} compares SIR/SAVE to LSIR/LSAVE in terms of sampling. Similar to the study in Section \ref{subsec:1d_ridge}, we use the same Monte Carlo samples for the slice-based and Lanczos-Stieltjes methods with various values of $R$ and $k$, respectively. The Lanczos-Stieltjes plots contain the best case from the slice-based algorithms for each value of $N$ for comparison. The new algorithms perform as well as the best case of the slice-based algorithms for all values of $k$ except for $k = 5$. For large values of $N$, the error associated with 5 Lanczos iterations levels off as the errors due to sampling become smaller than errors due to truncating the pseudospectral expansion. However, by increasing the number of Lanczos iterations, we improve the accuracy to the best case slice-based results.

%% file: sec7-conclusion.tex
\section{Conclusion}
\label{sec:conclusion}

We propose alternative approaches to the sliced inverse regression (SIR) and sliced average variance estimation (SAVE) algorithms for approximating $\CIR$ and $\CAVE$. The traditional methods approximate these matrices by applying a sliced partitioning over the range of output values. In the context of deterministic functions, this slice-based approach can be interpreted as a Riemann sum approximation of the integrals in \eqref{eq:CIR_CAVE}. The proposed algorithms use tools from classical numerical analysis, including orthonormal polynomials and Gauss-Christoffel quadrature, to produce high-order approximations of $\CIR$ and $\CAVE$. We refer to the new algorithms as Lanczos-Stieltjes inverse regression (LSIR) and Lanczos-Stieltjes average variance estimation (LSAVE).

We numerically study the various aspects of the Lanczos-Stieltjes algorithms on three test problems. We first examine the convergence of the approximate quadrature and orthonormal polynomial components resulting from the discrete approximation of the Stieltjes procedure by the Lanczos algorithm. This study highlights the interplay between the number of quadrature nodes over the input space and the number of Lanczos iterations. More Lanczos iterations correspond to higher degree polynomials, which require more quadrature nodes to accurately estimate. Poor approximations of these polynomials lead to poor approximations of $\CIR$ and $\CAVE$. We then compared the Lanczos-Stieltjes approximations of $\CIR$ and $\CAVE$ to their slice-based counterparts. These numerical studies emphasize a key characteristic of the Lanczos-Stieltjes methodology. Due to the composite structure of $\CIR$ and $\CAVE$, both the slicing approach and Lanczos-Stieltjes contain two levels of approximation: (i) sampling over the input space $\sX$ and (ii) approximation of the output space $\sF$. The Riemann sum approximation on $\sF$ results in a constant trade-off between (i) and (ii) in terms of which approximation is the dominant source of numerical error for various choices of $N$ and $R$. The Lanczos-Stieltjes approach significantly reduces the errors due to approximation over $\sF$, placing the burden of accuracy on the approximation over $\sX$. This enables Gauss-Christoffel quadrature over $\sX$ to produce the expected exponential convergence. Furthermore, when tensor product quadrature rules are infeasible, the Lanczos-Stieltjes approach allows Monte Carlo sampling to perform as expected without significant dependence on the approximation over $\sF$.